\documentclass[11pt]{article}

\usepackage{amsmath, amssymb, amscd}

\def\eqref#1{(\ref{#1})}
\newcommand{\goth}{\mathfrak}

\newcommand{\arrow}{{\:\longrightarrow\:}}
\newcommand{\Z}{{\Bbb Z}}
\newcommand{\C}{{\Bbb C}}
\newcommand{\R}{{\Bbb R}}
\newcommand{\Q}{{\Bbb Q}}

\newcommand{\6}{\partial}
\def\1{\sqrt{-1}\:}
\newcommand{\restrict}[1]{{\left|_{{\phantom{|}\!\!}_{#1}}\right.}}
\newcommand{\cntrct}                
{\hspace{2pt}\raisebox{1pt}{\text{$\lrcorner$}}\hspace{2pt}}

\def\Bbb#1{\mathbb #1}

\renewcommand{\bar}{\overline}
\renewcommand{\phi}{\varphi}
\renewcommand{\epsilon}{\varepsilon}

\renewcommand{\leq}{\leqslant}

\newcommand{\even}{{\rm even}}

\newcommand{\odd}{{\rm odd}}

\newcommand{\im}{\operatorname{im}}
\newcommand{\End}{\operatorname{End}}

\newcommand{\Hol}{\operatorname{{\cal H}ol}}

\newcommand{\Id}{\operatorname{Id}}

\newcommand{\Vol}{\operatorname{Vol}}

\renewcommand{\Re}{\operatorname{Re}}

\newcommand{\comment}[1]{{}}

\def\blacksquare{\hbox{\vrule width 4pt height 4pt depth 0pt}}
\def\endproof{\blacksquare}

\makeatletter

\newcommand{\ps@verbit}{%
  \renewcommand{\@oddhead}{%
          \scriptsize
          {Rational homotopy and parallel forms}
          \hfil\tiny {M. Verbitsky, February 23, 2005}}
  \renewcommand{\@evenhead}{\@oddhead}
  \renewcommand{\@oddfoot}{\hfil\thepage\hfil}
  \renewcommand{\@evenfoot}{\@oddfoot}}
 
\newcounter{Mycounter}[section]
\newcounter{lemma}[section]
\renewcommand{\thelemma}{\noindent{Lemma \thesection.\arabic{lemma}}}
\newcommand{\lemma}{%
     \setcounter{lemma}{\value{Mycounter}}
     \refstepcounter{lemma}
     \stepcounter{Mycounter}
     {\bf \thelemma:\ }}

\newcounter{claim}[section]
\renewcommand{\theclaim}{\noindent{Claim \thesection.\arabic{claim}}}
\newcommand{\claim}{%
     \setcounter{claim}{\value{Mycounter}}
     \refstepcounter{claim}
     \stepcounter{Mycounter}
     {\bf \theclaim:\ }}

\newcounter{sublemma}[section]
\newcounter{corollary}[section]
\renewcommand{\thecorollary}{\noindent{Corollary \thesection.\arabic{corollary}}}
\newcommand{\corollary}{%
     \setcounter{corollary}{\value{Mycounter}}
     \refstepcounter{corollary}
     \stepcounter{Mycounter}
     {\bf \thecorollary:\ }}

\newcounter{theorem}[section]
\renewcommand{\thetheorem}{\noindent{Theorem \thesection.\arabic{theorem}}}
\newcommand{\theorem}{%
     \setcounter{theorem}{\value{Mycounter}}
     \refstepcounter{theorem}
     \stepcounter{Mycounter}
     {\bf \thetheorem:\ }}

\newcounter{conjecture}[section]
\newcounter{proposition}[section]
\renewcommand{\theproposition}
       {\noindent{Proposition \thesection.\arabic{proposition}}}
\newcommand{\proposition}{%
     \setcounter{proposition}{\value{Mycounter}}
     \refstepcounter{proposition}
     \stepcounter{Mycounter}
     {\bf \theproposition:\ }}

\newcounter{definition}[section]
\renewcommand{\thedefinition}
       {\noindent{Definition~\thesection.\arabic{definition}}}
\newcommand{\definition}{%
     \setcounter{definition}{\value{Mycounter}}
     \refstepcounter{definition}
     \stepcounter{Mycounter}
     {\bf \thedefinition:\ }}

\newcounter{example}[section]
\newcounter{remark}[section]
\renewcommand{\theremark}{\noindent{Remark \thesection.\arabic{remark}}}
\newcommand{\remark}{%
     \setcounter{remark}{\value{Mycounter}}
     \refstepcounter{remark}
     \stepcounter{Mycounter}
     {\bf \theremark:\ }}

\newcounter{problem}[section]
\newcounter{question}[section]
\@addtoreset{equation}{section}
\@addtoreset{footnote}{section}
\makeatother

\begin{document}

\begin{center}
{\LARGE\bf
Manifolds with parallel differential forms \\[3mm]
and K\"ahler identities for $G_2$-manifolds
}
\\[4mm]
Misha Verbitsky\footnote{Misha Verbitsky is 
an EPSRC advanced fellow 
supported by CRDF grant RM1-2354-MO02 
and EPSRC grant  GR/R77773/01}
\\[4mm]

{\tt verbit@maths.gla.ac.uk, \ \  verbit@mccme.ru}
\end{center}

{\small 
\hspace{0.15\linewidth}
\begin{minipage}[t]{0.7\linewidth}
{\bf Abstract} \\
Let $M$ be a compact Riemannian manifold
equipped with a parallel differential form $\omega$.
We prove a version of K\"ahler identities in this
setting. This is used to
show that the de Rham algebra of $M$ is weakly  equivalent
to its subquotient $(H^*_c(M), d)$, called {\bf the pseudocohomology}
of $M$. When $M$ is compact and K\"ahler and $\omega$ is
its  K\"ahler form, $(H^*_c(M), d)$ is isomorphic to the 
cohomology algebra of $M$.
This gives another proof of homotopy
formality for K\"ahler manifolds, originally shown
by Deligne, Griffiths, Morgan and Sullivan.
We compute $H^*_c(M)$ for a compact $G_2$-manifold,
showing that $H^i_c(M)\cong H^i(M)$ unless 
$i=3,4$. For $i=3,4$, we compute
$H^*_c(M)$ explicitly in terms of the
first order differential operator
$*d:\; \Lambda^3(M)\arrow \Lambda^3(M)$.
\end{minipage}
}

{
\small
\tableofcontents
}


\section{Introduction}
\label{_Intro_Section_}


\subsection{Holonomy groups in Riemannian geometry}

Let $M$ be a Riemannian manifold equipped with a 
differential form $\omega$. This form is called
{\bf parallel} if $\omega$ is preserved by the
Levi-Civita connection: $\nabla\omega=0$.
This identity gives a powerful restriction
on the holonomy group $\Hol(M)$.

The structure of $\Hol(M)$ and its relation
to geometry of a manifold is one of the main subjects
of Riemannian geometry of last 50 years. This
group is compact, hence reductive,
and acts, in a natural way, on the tangent space $TM$.
When $M$ is complete,
Georges de Rham proved that unless this representation is irreducible,
$M$ has a finite covering, which is a product of Riemannian
manifolds of smaller dimension (\cite{_de_Rham_}). 
Irreducible holonomies were
classified by M. Berger (\cite{_Berger:list_}),
who gave a complete list of all irreducible 
holonomies which can occur on non-symmetric spaces.
This list is quite short:

\hfill

\begin{tabular}{|l|l|}
\hline
Holonomy  & Geometry\\[1mm]
\hline
$SO(n)$ acting on $\R^n$ & Riemannian manifolds\\[1mm]
\hline
$U(n)$ acting on $\R^{2n}$ & K\"ahler manifolds\\[1mm]
\hline
$SU(n)$ acting on $\R^{2n}$, $n>2$ & Calabi-Yau manifolds\\[1mm]
\hline
$Sp(n)$ acting on $\R^{4n}$ & hyperk\"ahler manifolds\\[1mm]
\hline
$Sp(n)\times Sp(1)/\{\pm 1\}$ & 
quaternionic-K\"ahler\\[1mm] acting on $\R^{4n}$, $n>1$ &  manifolds\\[1mm]
\hline
$G_2$ acting on $\R^7$ & $G_2$-manifolds \\[1mm]
\hline
$Spin(7)$ acting on $\R^8$ & $Spin(7)$-manifolds\\[1mm]
\hline
\end{tabular}

\hfill

\hfill

Berger's list also included $Spin(9)$ acting on $\R^{16}$,
but D. Alekseevsky later observed that this case is impossible
(\cite{_Alekseevskij:Spin(9)_}), unless $M$ is symmetric. 
If an irreducible manifold
$M$ has a parallel differential form, its holonomy is restricted,
as $SO(n)$ has no invariants in $\Lambda^i(TM)$, $0<i<n$.
Then $M$ is locally a product of symmetric spaces and
manifolds with holonomy $U(n)$, $SU(n)$, $Sp(n)$, etc.

\hfill

In K\"ahler geometry (holonomy $U(n)$) the parallel
forms are the K\"ahler form and its powers. Studying 
the corresponding algebraic structures, the algebraic geometers
amassed an amazing wealth of topological and geometric
information. In this paper we try to generalize some of
these results to other manifolds with a parallel form,
especially the $G_2$-manifolds. The results thus obtained
can be summarized as ``K\"ahler identities for
$G_2$-manifolds''.

\subsection{$G_2$-manifolds in mathematics and physics}

The theory of $G_2$-manifolds is one of the places where
mathematics and physics interact most intensely. For many years
after Berger's groundbreaking results, this subject was dormant;
after Alekseevsky showed that $Spin(9)$ cannot be realized in 
holonomy, there were doubts whether the other two exceptional
entries in Berger's list ($G_2$ and $Spin(7)$) can be realized.

Only in 1980-ies were manifolds with holonomy $G_2$
constructed. R. Bryant (\cite{_Bryant:holonomy_constru_})
found local examples, and then R. Bryant and S. Salamon
found complete manifolds with holonomy $G_2$
(\cite{_Bryant_Salamon_}). The compact examples of
holonomy $G_2$ and $Spin(7)$-manifolds were produced
by D. Joyce (\cite{_Joyce_G2_}, \cite{_Joyce_Book_}),
using difficult (but beautiful and quite powerful)
arguments from analysis and PDE theory. Since then,
the $G_2$-manifolds became a central subject of study
in some areas of string physics, and especially
in M-theory. The mathematical study of $G_2$-geometry
was less intensive, but still quite fruitful. 
Important results were obtained in 
gauge theory on $G_2$-manifold (the study of
Donaldson-Thomas bundles): \cite{_Donaldson_Thomas_},
\cite{_Tian:Calibrated_}, \cite{_Tian_Tao_}. 
A. Kovalev found many new examples of 
$G_2$-manifolds, using a refined version
of Joyce's engine (\cite{_Kovalev_}). 
N. Hitchin constructed a geometric flow
(\cite{_Hitchin:3-forms_}, \cite{_Hitchin:stable_}), which turned
out to be extremely important in string physics
(physicists call this flow {\bf Hitchin's flow}).
Hitchin's flow acts on the space of all ``stable''
(non-degenerate and positive) 3-forms on a 7-manifold.
It is fixed precisely on the 3-forms corresponding
to the connections with holonomy in $G_2$. 
In \cite{_Dijkgraaf_etc:M-theory_}, a
unified theory of gravity is introduced, based 
in part on Hitchin's flow. From the special cases
of topological M-theory one can deduce 4-dimensional
loop gravity, and 6-dimensional A- and B-models
in string theory. 

In string theory, $G_2$-manifolds are expected to play the
same role as Calabi-Yau manifolds in the usual
A- and B-model of type-II string theories. 
These two forms of string theory
both use Calabi-Yau manifolds, in a different 
fashion. Duality between these theories leads
to duality between Calabi-Yau manifolds, and then 
to far-reaching consequences, which were 
studied in mathematics and physics, under the
name of Mirror Symmetry. During the last 
20 years, the Mirror Symmetry became
one of the central topics of modern algebraic 
geometry.

There are two important ingredients in 
Mirror Symmetry (in Strominger-Yau-Zaslow form) -
one counts holomorphic curves on one Calabi-Yau
manifolds, and the special Lagrangian cycles on
its mirror dual. Using $G_2$-geometry, these 
two kinds of objects (holomorphic curves
and special Lagrangian cycles) are 
transformed into the same kind of objects,
called {\bf associative cycles} on a $G_2$-manifold.
This is done as follows.

A $G_2$-structure on a 7-manifold is given by
a 3-form (see Subsection \ref{_G_2_Subsection_}).
Consider a Calabi-Yau manifold $X$, $\dim M=3$, with 
non-degenerate holomorphic 3-form $\Omega$,
and Kaehler form $\omega$. Let $M:= X \times S^1$,
and let $dt$ denote the unit cotangent form
of $S^1$ lifted to $M$. Consider a 3-form
$\omega \wedge dt + \Re \Omega$ on $M$.
This form is obviously closed. 
It is easy to check that it defines a 
parallel $G_2$-structure on $M$.
This way one can convert problems from
Calabi-Yau geometry to problems
in $G_2$-geometry. 

A 3-form $\phi$ on a manifold 
$M$ gives an anti-symmetric
map
\[ \phi^\sharp:\; TM \otimes TM \arrow \Lambda^1(M),
\]
$x, y \arrow \phi(x, y , \cdot)$.
Using the Riemannian structure, we identify
$TM$ and $\Lambda^1(M)$. Then
$\phi^\sharp$ leads to a skew-symmetric 
vector product $V:\; TM \otimes TM \arrow TM$.
An {\bf associative cycle} on a $G_2$-manifold 
is a 3-dimensional submanifold $Z$ such that 
$TZ$ is closed under this vector product. 
Associative submanifolds are studied within
the general framework of calibrated geometries 
(see \cite{_Harvey_Lawson:Calibrated_}).

Given a Calabi-Yau threefold $X$,
consider $M= X\times S^1$ with a $G_2$-structure
defined above. Let $Z\subset X$ be a 3-dimensional
submanifold. It is easy to check
that $Z$ is special Lagrangian if and only if
$Z\times \{t\}$ is associative in $M$.
Also, given a 2-cycle $C$ on $X$,
$C\times S^1$ is associative in $M$ if
and only if $C$ is a holomorphic curve.
This way, the instanton objects in mirror dual
theories (holomorphic curves and SpLag cycles)
can be studied uniformly after passing to $G_2$-manifold.
It was suggested that this correspondence indicates
some  form of string duality (\cite{_Leung:TQFT_}, \cite{_Salur:Santillan_}).

However, the main physical motivation for the study
of $G_2$-manifolds comes from M-theory; we direct the
reader to the excellent survey \cite{_Acharya_Gukov_} for
details and further reading. M-theory
is a theory which is expected (if developed) to produce
a unification of GUT (the Grand Unified Theory
of strong, weak and electro-magnetic forces)
with gravity, via supersymmetry. In this
approach, string theories arise as approximations
of M-theory. In most applications related to
M-theory, a $G_2$-manifold is deformed to a 
compact $G_2$-variety with isolated singularities.
One local construction of conical singularities
of this type is based on Bryant-Salamon examples
of complete $G_2$-manifolds (see \cite{_Bryant_Salamon_}).
In this approach, the study of conical singularities 
is essentially reduced to the 4-dimensional geometry.

An explicit mathematical study of these singular examples
and their connection to physics and theory
of Einstein manifolds is found in 
\cite{_Atiyah_Witten_}. Also, Hitchin's
flow can be used to produce many such
examples in a uniform way (see 
\cite{_Gukov_Yau_Zaslow_})

\subsection{Structure operator on manifolds with
parallel differential form}

Much study in K\"ahler geometry is based on the interplay
between the
de Rham differential and the twisted de Rham differential
$d^c := - I \circ d \circ I$. We construct a similar
operator $d_c$ for any manifold with a parallel differential
form. This operator no longer satisfies $d_c^2=0$;
however, it satisfies many properties expected
from the twisted de Rham differential in K\"ahler
geometry. Most importantly, a version of $dd_c$-lemma
is true in this setting (\ref{_d_d_c_Proposition_}).

Just as in the usual case, this may lead
to results in rational homotopy theory (see Subsection 
\ref{_formality_G_2_Subsection_} in the present introduction).

To simplify the exposition, we restrict ourselves 
presently to Riemannian manifolds $(M, \omega)$
with a parallel 3-form. These include 
Riemannian 3-manifolds, Calabi-Yau threefolds
and $G_2$-manifolds. Just like it happens
in 3-dimensional case, such a 3-form defines a
skew-symmetric cross-product on $\Lambda^1(M)$:
\[
x, y \stackrel\Psi \arrow \omega(x^\sharp, y^\sharp, \cdot)
\]
($(\cdot)^\sharp$ denotes taking the dual with respect to
the metric).
Consider the operator on differential forms
\begin{align*}
\xi_{i_1}\wedge & \xi_{i_2} \wedge ... \wedge \xi_{i_k}\\
\arrow & \sum_{1\leq a <b \leq k} 
       (-1)^{a+b-1}\Psi(\xi_{i_a},\xi_{i_b})
        \wedge \xi_{i_1}\wedge \xi_{i_2} \wedge ... \wedge 
       \hat\xi_{i_a}\wedge ...\wedge 
       \hat\xi_{i_b}\wedge ... \wedge \xi_{i_k}
\end{align*}
where $\xi_i$ is an orthonormal frame in $\Lambda^1(M)$.
Denote by \[ C:\; \Lambda^i(M) \arrow \Lambda^{i+1}(M)\]
the dual operator (to identify $\Lambda^i(M)$
with its dual, we use the natural metric
on $\Lambda^i(M)$ induced from the Riemannian
structure on $M$). Then $C$ is called {\bf
the structure operator on $(M,\omega)$}.

In Section \ref{_para_form_Section_}
we give another definition of $C$,
which works for an arbitrary parallel 
$i$-form $\omega$. It is not difficult to check
that this definition is compatible to 
the one given above. When $(M, \omega)$ is K\"ahler,
$C$ becomes the complex structure operator
on $M$, and the identities we prove
in general case become the usual 
K\"ahler identities.

Denote by $d_c$ the anticommutator $\{C, d\}= dC + Cd$.
We show that $d_c$ commutes with $d$, $d^*$, and
satisfies the following version of $dd_c$-lemma

\hfill

\proposition\label{_d_d_c_Proposition_}
Consider a compact Riemannian manifold
equipped with a parallel differential form.
Let $\eta$ be a differential $k$-form satisfying
$d\eta = d_c \eta =0$. 
Assume, moreover, that $\eta$ is $d_c$-exact:
$\eta = d_c \xi$. Then $\eta = d d_c \xi'$, 
for some differential form $\xi$.

\hfill

{\bf Proof:} Follows immediately from
\ref{_H_c^*_quasi-iso_Proposition_} (see \ref{_ddc_proof_Remark_}).
 \endproof

\hfill

\remark 
The operator $d_c$ satisfies the Leibniz identity:
\[ d_c(a \wedge b) = d_c(a) \wedge b + 
   (-1)^{\tilde a \tilde d_c} a \wedge d_c (b),
\]
where $\tilde a$, $\tilde b$ denotes parity of a form.
However, $d_c^2\neq 0$. Also, the $d d_c$-lemma is less
strong than the usual $d d^c$-lemma: given a $d$-exact,
$d, d_c$-closed form $\eta$, we cannot show that
$\eta = d d_c \xi'$ (though this could be true
in the case of $G_2$-manifolds).

\subsection{Donaldson-Thomas bundles}

The twisted de Rham operator has many uses in
$G_2$-geometry. In many ways, $d_c$ defines the same kind
of structures as known in algebraic geometry from the
study of the holomorphic structure operator 
$\bar \6 = \frac{d - \1 d^c}2$.

Let $M$ be a $G_2$-manifold. The $G_2$-action gives
a decomposition
\[ \Lambda^2 (M) = \Lambda^2_7(M)\oplus \Lambda^2_{14}(M)
\]
onto a sum of irreducible representations of $G_2$.

\hfill

\definition
\cite{_Donaldson_Thomas_}
Let $(B, \nabla)$ be a vector bundle with connection
on a $G_2$-manifold $M$, and 
$\Theta \in \Lambda^2(M)\otimes \End(B)$ 
its curvature. Then $(B, \nabla)$
is called {\bf a Donaldson-Thomas bundle}
if $\Theta$ lies inside 
$\Lambda^2_{14}(M)\otimes \End(B)$.

\hfill

This is a natural generalization of the Hermitian-Einstein
condition, known from algebraic geometry.
In fact, when $M$ is constructed from a Calabi-Yau
threefold $W$, $M = W \times S^1$, the Donaldson-Thomas
bundles can be obtained as a pullback of Hermitian-Einstein
bundles from $W$ to $M$. Also, the Donaldson-Thomas
condition implies that the functional
\[ (B, \nabla) \arrow \int_M ||\Theta||^2\Vol(M) \]
has an absolute minimum at $(B, \nabla)$. In other words, 
Donaldson-Thomas bundles are always instantons.

Geometry of Donaldson-Thomas bundles is much
studied in connection with physics and algebraic
geometry, see e.g. \cite{_Leung:TQFT_}, 
\cite{_Lee_Leung_}.  

Given a Hermitian vector bundle $(B, \nabla)$
on a K\"ahler manifold, the holomorphic condition can be written
as $\Theta \in \Lambda^{1,1}(M)\otimes \End(B)$.
This equation can be rewritten as $\{\nabla, \nabla^c\}=0$,
where 
$\nabla^c = - I \circ \nabla \circ I = [ W_I,\nabla]$
($W_I$ denotes the K\"ahler-Weil operator, acting on
$\Lambda^{p,q}(M)$ as $\1 (p-q)$). In $G_2$-geometry
the role of $W_I$ is played by the structure 
operator $C$.

The Donaldson-Thomas bundles can be interpreted in terms
of a structure operator, repeating the above desctiption for
holomorphic bundles verbatim. 

\hfill

\proposition
Let $M$ be a $G_2$-manifold, 
$C:\; \Lambda^i(M) \arrow\Lambda^{i+1}(M)$
the structure operator, and $(B, \nabla)$ a vector bundle
with connection,
\[ 
\nabla:\; B\otimes \Lambda^i(M) \arrow B\otimes \Lambda^{i+1}(M).
\]
Consided an operator $\nabla^c\;:= \{ C, \nabla\}$,
\[ 
\nabla^c:\; B\otimes \Lambda^i(M) \arrow B\otimes \Lambda^{i+2}(M).
\] 
Then $(B, \nabla)$ is a Donaldson-Thomas bundle if and
only if $\nabla$, $\nabla^c$ commute.

\hfill

{\bf Proof:} Using graded Jacobi idenity, we obtain
\[ [ C, \nabla^2] = \frac 1 2 [C, \{\nabla, \nabla\}]=
   [\nabla, [C, \nabla]]= [\nabla, \nabla^c].
\]
However, $[ C, \nabla^2]= C(\Theta)$, where $\Theta$
is the curvature form. In
\ref{_C_isomo_on_Lambda_7_Proposition_} 
we show that $\ker C\restrict{\Lambda^2(M)}$
is exactly $\Lambda^2_{14}(M)$, hence 
\[ C(\Theta)=0 \ \ \Leftrightarrow\ \ 
  \Theta\in \Lambda^2_{14}(M)\otimes\End(B).
\]
\endproof

\subsection{Localization functor and rational homotopy}
\label{_rational_homo_Subsection_}

The homotopy formality for K\"ahler manifold,
observed by Deligne,  Griffiths,  Morgan, 
Sullivan (\cite{_DGMS:Formality_}), is one of the
deepest and most powerful results of K\"ahler geometry.
Since \cite{_DGMS:Formality_} appeared, there was
a whole cornucopia of research dedicated to this theme.
Formality was used to study the deformations and 
moduli spaces (see e.g. \cite{Goldman_Millson_},
\cite{_Barannikov_Kontsevich_}, 
\cite{_Verbitsky:Hyperholo_bundles_}), 
in Mirror Symmetry and topology. 
The reason for all these equations lies
in the so-called {\bf Master equation}
(also known as {\bf the Maurer-Cartan equation})
\[ d \gamma = -\frac 1 2 [\gamma, \gamma]. \]
in a differential graded (DG-) Lie algebra,
which is responsible for deformation theory for
most objects in algebraic geometry.
Solutions of this equation (up to
a relevant equivalence) are homotopy
invariants of the DG-Lie algebra
(\cite{_Barannikov_Kontsevich_}).

This equation can be solved recursively,
if the relevant Massey products vanish
(in fact, Massey products can be defined
as obstructions to finding solutions of 
Maurer-Cartan equation - see e.g.
\cite{_Babenko_Taimanov_}). The homotopy
formality implies vanishing of Massey products,
providing a way to solve the 
Maurer-Cartan equation in 
various contexts. 

In the proof of homotopy formality 
for K\"ahler manifolds (\cite{_DGMS:Formality_}), the key
argument hinges on $dd^c$-lemma; one should expect that
the $G_2$-version of $dd^c$-lemma (\ref{_d_d_c_Proposition_})
will give us information about rational homotopy of
$G_2$-manifolds.

The topological utility of rational homotopy is based on
the Quillen-Sullivan localization construction,
\cite{_Quillen_}, \cite{_Sullivan:geom_top_}.
The $\Q$-localization functor in homotopy category
maps a simply connected cellular 
space $X$ to a space $X_\Q=Loc_\Q(X)$ with 
$H^i(X_\Q, \Z) \cong H^i(X, \Z)\otimes \Q$ and
$\pi_i(X_\Q) \cong \pi_i(X)\otimes \Q$.
The spaces which are homotopy equivalent
to their localization are called {\bf $\Q$-local}.
We have $Loc_\Q(X)\cong Loc_\Q(Loc_\Q(X))$;
in other words, all spaces of form
$Loc_\Q(X)$ are $\Q$-local.

Given a cellular space, one could construct
its de Rham complex, using piecewise smooth 
diferential forms. This construction maps
homotopy equivalent spaces to weakly equivalent
differential graded (DG-) 
algebras (see \ref{_formal_DG_Definition_}).
We obtain a functor 
$DR:\; \operatorname{Hot} \arrow \operatorname{DG-Alg}$
of the corresponding categories.
Moreover, this functor commutes with
localization, and gives  an equivalence
of homotopy category of $\Q$-local simply 
connected spaces and the category $\operatorname{DG-Alg}$
of DG-algebras. This reduces the study of rational homotopies
(homotopies of $\Q$-local spaces) to the
study of  DG-algebras. 

The localization construction (which is defined
in many other contexts, see \cite{_Dwyer:Localization_}) 
is one of the key ideas of modern
algebraic topology. Sullivan needed localization
in order to prove the Adams' conjecture, and Quillen
used localization to give the definition of algebraic K-theory.
Since then, many other uses of the same construction
were found; including Voevodsky's celebrated motivic homotopy theory.

Two DG-algebras are called {\bf quasi-isomorphic}
if there exists a quasi-isomorphism (morphism, inducing
isomorphism on cohomology) from one to another. 
The equivalence relation generated by 
quasi-isomorphism is called {\bf weak equivalence}
of DG-algebras (\ref{_formal_DG_Definition_}).

Rational homotopy is a study of DG-algebras,
up to weak equivalence.

A DG-algebra $(A^*, d)$ is called {\bf homotopy formal}
if it is weakly equivalent to its cohomology algebra
$(H^*(A), 0)$. A simply connected topological
space is called {\bf formal} if its de Rham algebra is formal.
The rational homotopies of formal spaces (in particular,
all rational homotopy groups) are determined by the
algebraic structure on cohomology.

Not all DG-algebras are formal; the best known obstruction
to formality is called {\bf the Massey product}
(see e.g. \cite{_Babenko_Taimanov_}). However,
there are more obstructions to formality
that just a Massey product.  S. Halperin and J. Stasheff 
(\cite{_Halperin_Stasheff_}) constructed explicitly 
a complete set of obstructions 
\[ \{ O_n, n =1, 2, 3, ...\} \]
to homotopy formality,  $O_n$ defined if all
$O_i$, $i<n$ vanish. 

Since homotopy formality of K\"ahler manifolds was established,
many people studied the influence of differential geometric structures
on rational homotopy. Much of this work
was focused on the study of rational homotopy of 
compact symplectic manifolds (there is a book 
\cite{_Tralle_Oprea_}, dedicated especially to this subject).
Using Deligne-Griffiths-Morgan-Sullivan formality theorem,
 one obtains all kinds of symplectic manifolds
admitting no K\"ahler structures. 

In physics,
$G_2$-manifolds appear as a generalization
of Calabi-Yau threefolds; formality is expected.

\subsection{Formality for $G_2$-manifolds}
\label{_formality_G_2_Subsection_}

Homotopy formality for $G_2$-manifolds was studied by Gil Cavalcanti
in his thesis (see \cite{_Cavalcanti_}). The $G_2$-structure gives certain
constraints on the cohomology ring of 
a manifold: the multiplication by the standartd 3-form
$\omega$ gives an isomorphism 
\[ H^2(M) \stackrel {\wedge \omega} \arrow H^5(M).
\]
and the 2-form 
\[
\eta \arrow \int_M \eta\wedge\eta \wedge \omega
\]
on $H^2(M, \R)$ must be positive definite.
Also, $H^1(M)=0$. Cavalcanti constructed
examples of non-formal 7-manifolds satisfying
these constraints. He also showed that for
$\dim H^2(M)\leq 2$, these constraints indeed imply formality.

The $G_2$-version of $dd_c$-lemma (\ref{_d_d_c_Proposition_})
should give information about rational homotopy, in the same way that
the usual $dd^c$-lemma leads to formality of K\"ahler manifolds.
Indeed, $(\ker d_c)$ is a subalgebra
of $\Lambda^*(M)$ which is weakly equivalent to
the de Rham algebra of $M$ (\ref{_ker_quasi-iso_Proposition_}),
and the quotient algebra 
\[ (H_c^*(M), d) \cong \frac {\ker d_c} {(\ker d_c) \cap (\im d_c)}
\]
is also weakly equivalent to $\Lambda^*(M)$.
We call $(H_c^*(M), d)$ {\bf the pseu\-do-\-co\-ho\-mo\-lo\-gy} of $M$
(\ref{_pseudoco_Definition_}). We don't call it cohomology,
because $d^2_c\neq 0$.

A form $\eta\in \Lambda^*(M)$ is called {\bf pseudo-harmonic}
if $\eta \in (\ker d_c) \cap(\ker^* d_c)$,
where $d^*_c$ is a Hermitian adjoint to $d_c$.
Just as happens for usual cohomology,
the space of pseudo-harmonic forms ${\cal H}_c^*(M)$ is isomorphic to
pseudo-cohomology:
\[ (H_c^*(M), d)\cong ({\cal H}_c^*(M), d)\]
(\ref{_pseudoha_pseudoco_Proposition_}).
All harmonic forms are also pseudo-harmonic. 
We consider an orthogonal decomposition
\[ 
{\cal H}_c^*(M)\cong {\cal H}^*(M)\oplus {\cal H}_c^*(M)_{>0},
\]
where ${\cal H}_c^*(M)_{>0}$ is the sum of all positive eigenspaces
of the Laplacian acting on ${\cal H}_c^*(M)$.
From the arguments given above, we immediately obtain the
following theorem.

\hfill

\theorem 
Let $M$ be a compact $G_2$-manifold, and ${\cal H}_c^*(M)_{>0}$
the sum of all positive eigenspaces
of the Laplacian acting on ${\cal H}_c^*(M)$.
Assume that ${\cal H}_c^*(M)_{>0}=0$.
Then $M$ is formal.

\hfill

{\bf Proof:} This is \ref{_homotopi_equiv_for_pseudo_Corollary_}.
\endproof

\hfill

We were unable to show that ${\cal H}^*_c(M)_{>0}=0$
for all $G_2$-manifolds.
However, this space was computed fairly  
explicitly, in terms of $G_2$-action on
differential forms.

\hfill

\proposition\label{_H^i_c(M)_>0_Proposition_}
Let $M$ be a compact $G_2$-manifold, and
${\cal H}^i_c(M)_{>0}=0$ the vector space defined above.
Then ${\cal H}^i_c(M)_{>0}=0$ unless $i=3$ or $4$.
The space ${\cal H}^3_c(M)_{>0}$ is generated
(over $\C$) by the solutions of the following equation
\begin{equation}\label{_d_alpha=mu_*_alpha_Equation_}
d\alpha = \mu *\alpha, \ \ \alpha \in \Lambda^3_{27}(M), 
\end{equation}
where $\mu \in \C$ is a non-zero number, and
$\Lambda^3_{27}(M)$ is the 27-dimensional 
irreducible component of $\Lambda^3(M)$ under
the $G_2$-action (see \eqref{_Lambda^*_splits_Equation_}). 
Similarly, ${\cal H}^4_c(M)_{>0}$ is generated by the solutions
of equation $d * \eta = \mu \eta$, $\eta\in \Lambda^4_{27}(M)$.

\hfill

{\bf Proof:} See \ref{_pseudoha_G_2_Theorem_}. 
\endproof 

\hfill

The formula \eqref{_d_alpha=mu_*_alpha_Equation_}
is suggestive of equations found in 
Hitchin's paper on hamiltionian flow,
\cite{_Hitchin:stable_}. One may hope
that a careful study of Hitchin's flow
in conjunction with \eqref{_d_alpha=mu_*_alpha_Equation_}
leads to some constraints on ${\cal H}^3_c(M)_{>0}$,
and, possibly, its vanishing, which leads to formality
of $M$. However, even now \ref{_H^i_c(M)_>0_Proposition_}
gives us some information about rational homotopy.

\hfill

\corollary\label{_diffe_on_H^2_c_Corollary_}
Let $M$ be a compact $G_2$-manifold, and 
$(H_c^*(M), d)$ its pseudocohomology DG-algebra. Then
$(H_c^*(M), d)$  is weakly equivalent to the de Rham 
algebra of $M$, and, moreover, $d\restrict{H_c^i(M)}=0$
unless $i=3$. 

\noindent
\endproof

\hfill

This result can be used to study the obstructions
$O_n$ to formality of $(H_c^*(M), d)$, defined in
\cite{_Halperin_Stasheff_} (see Subsection 
\ref{_rational_homo_Subsection_}). 
It turns out that only the first obstruction
$O_1$ is relevant for rational homotopy,
and if it vanishes, $O_i$, $i>0$ also vanish, 
and  the DG-algebra $(H_c^*(M), d)$
and $(\Lambda^*(M), d)$ is formal.
However, the same result can be obtained
from Gil Cavalcanti's work, for all
simply connected 7-manifolds.

In 1970-is, T.~J.~Miller
showed that all simply connected orientable compact manifolds
of dimensions up to 6 are formal (\cite{_Miller:formality_6_7_}).
Moreover, Miller has shown that all $(k-1)$-connected
orientable compact manifolds of dimension up to $2k+2$
are formal. His arguments were simplified and generalized by 
M. Fernandez and V. Munoz (\cite{_Fernandez_Munoz_}),
who defined a notion of {\bf $k$-formal manifold}, and
shown that any orientable $k$-formal compact manifold
of dimension up to $2k+2$ is formal. They applied this
theorem to obtain results about formality of compact 
symplectic  manifolds.

G. Cavalcanti
(\cite{_Cavalcanti_}) studied 7-manifolds using 
the same conceptual framework, obtaining essentially 
(but in completely different terms) that obstructions 
to 3-formality for simply connected 7-manifolds
can be reduced to vanishing of the first 
obstruction of Halperin-Stasheff. 

It is unclear whether additional topological
information might be gleaned from 
\ref{_diffe_on_H^2_c_Corollary_}. It may
possibly happen that for any 7-manifold
(compact and oriented) its de Rham algebra
is weakly equivalent to an algebra with
non-degenerate Poincare pairing and
a differential which vanishes 
in all dimensions except $i=3$.
In this case we don't obtain much
topological information from
\ref{_diffe_on_H^2_c_Corollary_}.


\section{Riemannian manifolds with a parallel differential
form}
\label{_para_form_Section_}


\subsection{Structure operator and twisted differential}

Let $M$ be a $C^\infty$-manifold. We denote the
smooth forms on $M$ by $\Lambda^*(M)$.
Given an odd or even form $\alpha\in\Lambda^*(M)$, 
we denote by $\tilde \alpha$ its parity,
which is equal to 0 for even forms, and 1 for odd forms.
An operator $f\in \End(\Lambda^*(M))$ 
preserving parity is called {\bf even}, and one exchanging
odd and even forms is {\bf odd}; $\tilde f$ is equal
0 for even forms and 1 for odd.
 
Given a $C^\infty$-linear map 
$\Lambda^1(M) \stackrel {p} \arrow\Lambda^\odd(M)$
or $\Lambda^1(M) \stackrel p \arrow\Lambda^\even(M)$,
$p$ can be uniquely extended to a $C^\infty$-linear
derivation $\rho$ on $\Lambda^*(M)$, using the rule
\[ 
\rho\restrict{\Lambda^1(M)}=p,\ \
\rho\restrict{\Lambda^0(M)}=0, \ \  
\rho(\alpha\wedge \beta) = \rho(\alpha)\wedge \beta + 
(-1)^{\tilde \rho\tilde\alpha} \alpha\wedge \rho(\beta).
\]
Then, $\rho$ is an even (odd) differentiation of the
graded commutative algebra $\Lambda^*(M)$.

\hfill

\definition
Let $M$ be a Riemannian manifold, 
and $\omega\in \Lambda^k(M)$ a differential form.
Consider an operator
$\underline C:\; \Lambda^1(M) \arrow \Lambda^{k-1}(M)$
mapping $\nu\in \Lambda^1(M)$ to
$\omega\cntrct\nu^\sharp$, where 
$\nu^\sharp$ is the vector field 
dual to $\nu$. Alternatively,
$\underline C(\nu)$ can be written as
$\underline C(\nu) = *(*\omega\wedge \nu)$.
The corresponding differentiation
\[ C:\; \Lambda^*(M) \arrow \Lambda^{*+k-2}(M)
\]
is called {\bf the structure operator of $(M,\omega)$}.
Parity of $C$ is equal to that of $\omega$.

\hfill

\remark
When $(M,I,g)$ is a K\"ahler manifold 
and $\omega$ is its K\"ahler form, $\underline C(\nu) = I(\nu)$, 
and $C$ is the standard K\"ahler-Weil operator, acting on $(p,q)$-forms
as a multiplication by $(p-q)\1$.

\hfill

\definition
Let $M$ be a Riemannian manifold, 
$\omega\in \Lambda^k(M)$ a differential form,
which is parallel with respect to the Levi-Civita
connection. Denote by $d_c$ the supercommutator
\[
\{d, C\}:= dC - (-1)^{\tilde C} Cd
\]
This operator is called {\bf the twisted de Rham
operator of $(M,\omega)$}. Being a graded commutator
of two graded differentiations, $d_c$ is also a 
graded differentiation of $\Lambda^*(M)$.

\remark
When $(M,I,g)$ is a K\"ahler manifold 
and $\omega$ is its K\"ahler form, $d_c$ is
 equal to the well-known twisted differential 
$d^c = I^{-1}\circ d\circ I$, 
$d^c = \frac{\6-\bar \6}{\1}$.
Of course, for a general form $\omega$, 
$d_c^2$ can be non-zero.

\hfill

\proposition\label{_d_c_via_omega_Proposition_}
Let $(M, \omega)$ be a Riemannian manifold 
equipped with a parallel form $\omega$, and
$L_\omega$ the operator $\eta \arrow \eta \wedge \omega$.
Then 
\[
d_c = \{L_\omega, d^*\},
\]
where $\{\cdot,\cdot\}$ denotes the supercommutator,
\[ \{L_\omega, d^*\} = 
    L_\omega d^*- (-1)^{\tilde \omega}d^* L_\omega,
\]
and $d^*= -*d*$ is the adjoint to $d$.

\hfill

{\bf Proof:}
Denote by $\nabla$ the Levi-Civita connection,
\[ \nabla:\; \Lambda^*(M)\arrow\Lambda^*(M)\otimes
   \Lambda^1(M).
\] 
Let $\eta\in \Lambda^i(M)$. 
Clearly, $d^*\eta$ is obtained  from 
$\nabla \eta\in \Lambda^i(M)\otimes \Lambda^1(M)$
by applying the isomorphism 
\[ \Lambda^i(M)\otimes \Lambda^1(M)\cong 
   \Lambda^i(M)\otimes TM
\] 
induced by the Riemannian structure 
and then plugging the $TM$-part into $\Lambda^i(M)$:
\begin{equation}\label{_d^*_via_nabla_Equation_}
d^*\eta = \cntrct(\nabla \eta).
\end{equation}
Since $\nabla\omega=0$, $\{L_\omega, d^*\}$
is equal to the composition
\begin{align*}
\Lambda^i(M)&\stackrel\nabla\arrow
\Lambda^i(M)\otimes \Lambda^1(M) \\
\  & \stackrel{C\otimes \Id}\arrow 
\Lambda^{i+k-2}(M)\otimes \Lambda^1(M)
\stackrel\wedge  \arrow\Lambda^{i+k-1}(M)
\end{align*}
(the last arrow is exterior multiplication).
Indeed, $L_\omega$ commutes with $\nabla$, and
therefore, by \eqref{_d^*_via_nabla_Equation_},
$\{L_\omega, d^*\}$ is written
as a composition of $\nabla$ and a
commutator of $C^\infty$-linear maps
$L_\omega$ and $\cntrct$, where
\[ \cntrct:\; \Lambda^{i+k}(M)\otimes\Lambda^1(M)\arrow \Lambda^{i+k-1}(M)\]
maps $\eta\otimes \nu$ to $\eta\cntrct\nu^\sharp$. 
However, by definition,
\[
\{L_\omega, \cntrct\}(\eta\otimes \nu) = C(\eta) \wedge \nu.
\]
This gives
\begin{equation}\label{_commu_omega_contrct_Equation_} 
\{L_\omega, d^*\}(\eta) = [L_\omega, \cntrct](\nabla\eta).
\end{equation}

Similarly, $[\nabla, C]=0$, hence $d_c$ is written as a
composition of $\nabla$ and a $C^\infty$-linear map
\begin{equation}\label{_Id_C_wedge_commu_Equatuion_}
C\otimes \Id\circ \wedge - \wedge \circ C:\; 
\Lambda^{i}(M)\otimes \Lambda^1(M)\arrow \Lambda^{i+k-1}(M),
\end{equation}
where $\wedge:\; \Lambda^{*}(M)\Lambda^1(M)\arrow \Lambda^{*+1}(M)$
denotes the exterior product.
Since $C$ is a differentiation, the operator
\eqref{_Id_C_wedge_commu_Equatuion_} is equal to
\[ \Id\otimes C\circ \wedge:\; 
\Lambda^{i}(M)\otimes \Lambda^1(M)\arrow \Lambda^{i+k-1}(M).
\]
This gives
\begin{equation}\label{_commu_C_Equation_} 
\{d, C\}(\eta) = \Id\otimes C\circ \wedge(\nabla\eta).
\end{equation}
However,  by definition of $C$, we have
$[L_\omega, \cntrct](\eta\otimes\nu)=\eta\wedge C(\nu)$,
hence the right hand sides of \eqref{_commu_C_Equation_} 
and \eqref{_commu_omega_contrct_Equation_} are equal. 
This proves \ref{_d_c_via_omega_Proposition_}. \endproof

\hfill

\remark
In the K\"ahler case, \ref{_d_c_via_omega_Proposition_}
becomes the following well-known K\"ahler identity:
$[L_\omega, d^*]=d^c$.

\subsection{Generalized K\"ahler identities and twisted Laplacian}

\proposition\label{_kahler_ide_Proposition_}
Let $M$ be a Riemannian manifold equipped with a
parallel differential $k$-form $\omega$, 
$d_c$ the twisted de Rham operator constructed above,
and $d^*_c$ its Hermitian adjoint. Then
\begin{description}

\item[(i)] The following supercommutators vanish:
\[ \{ d, d_c\} =0, \ \ \{ d, d_c^*\} =0, \ \ 
   \{ d^*, d_c\} =0, \ \  \{ d^*, d_c^*\} =0, 
\]
\item[(ii)] The Laplacian $\Delta=\{d, d^*\}$
commutes with $L_\omega:\; \eta \arrow \omega\wedge \eta$
and its Hermitian adjoint operator, denoted
as $\Lambda_{\omega}:\; \Lambda^i(M)\arrow\Lambda^{i-k}(M)$.

\item[(iii)] Denote the supercommutator of $d_c, d^*_c$ by
$\Delta_c$. By definition, 
$\Delta_c = d_c d^*_c + d^*_c d_c$
when $k$ is even, and $\Delta_c = d_c d^*_c - d^*_c d_c$
when $k$ is odd. Then 
\[ 
 \Delta_c = (-1)^{\tilde \omega}\{ d^*, [H_\omega, d]\},
\]
where $H_\omega = \{ L_\omega, \Lambda_\omega\}$.
\end{description}

\hfill

{\bf Proof:} We use the following basic lemma

\hfill

{\bf Basic Lemma:} Let $\delta$ be an odd element in a
graded Lie superalgebra $A$  satisfying $\{\delta, \delta\}=0$.
Then $\{\delta, \{\delta, x\}\}=0$ for all $x\in A$, 
assuming that the base field is not of 
characteristic 2.

\hfill

{\bf Proof:} Using the graded Jacobi identity, we obtain
\[ \{\delta, \{\delta, x\}\} = - \{\delta, \{\delta, x\}\} + \{\{\delta, \delta\}, x\}. 
\]
This gives $2\{\delta, \{\delta, x\}\}=0$.
\endproof

\hfill

Now, $\{d, d_c\} = \{d, \{d, C\}\}=0$ (by the Basic Lemma),
and $\{d^*, d_c\} = \{d^*, \{d^*, L_{\omega}\}\}=0$
(by Basic Lemma and \ref{_d_c_via_omega_Proposition_}).
Taking Hermitian adjoints of these identities, 
we obtain the other two equations of 
\ref{_kahler_ide_Proposition_} (i).
\ref{_kahler_ide_Proposition_} (i) is proven.

Now, the graded Jacobi identity implies
\begin{equation}\label{_Laplacian_commu_L_omega_Equation_}
[ L_\omega, \Delta] =\{ L_\omega, \{d, d^*\}\} =
(-1)^{\tilde \omega}\{ d, \{L_\omega, d^*\}\}.
\end{equation}
(we use $\{L_\omega, d\}=0$ as $\omega$ is closed).
This gives
\[ 
  [L_\omega, \Delta] = (-1)^{\tilde \omega}\{ d, d_c\}=0,
\]
as \ref{_kahler_ide_Proposition_} (i) implies. 
Taking Hermitian adjoint, we also obtain
$[\Lambda_\omega, \Delta]=0$. We proved
\ref{_kahler_ide_Proposition_} (ii).

Finally, \ref{_kahler_ide_Proposition_} (iii)
is proven as follows:
\begin{equation} \label{_Delta_c_commu_Equation_}
\{\{L_\omega, d^*\}, \{\Lambda_\omega, d\}\}=
\{\{L_\omega, \{d^*, d_c^*\}\}+ 
(-1)^{\tilde \omega} \{d^*, \{L_\omega, \{\Lambda_\omega, d\}\}\}
\end{equation}
by graded Jacobi identity. Also,
\begin{equation}\label{_L_Lambda_d_Equation_}
   \{ L_\omega,\{\Lambda_\omega, d\}\} = \{ H_\omega, d\}
       +(-1)^{\tilde \omega}\{\Lambda_\omega, \{L_\omega, d\}\}.
\end{equation}
However, $\{L_\omega, d\}=0$ as $\omega$ is
closed. Comparing \eqref{_L_Lambda_d_Equation_} and
\eqref{_Delta_c_commu_Equation_}, we obtain
\[
\Delta_c = (-1)^{\tilde\omega} \{d^*, \{ H_\omega, d\}\}.
\]
We proved \ref{_kahler_ide_Proposition_} (iii).
\endproof

\hfill

\remark 
When $(M, \omega)$ is a K\"ahler manifild, 
\ref{_kahler_ide_Proposition_} (i) gives the standard
commutation relations between $d$, $d^c$, $d^*$, $(d^c)^*$,
\ref{_kahler_ide_Proposition_} (ii) is well known, and 
\ref{_kahler_ide_Proposition_} (iii) gives 
\[
 \{d^c,(d^c)^*\}= \Delta_c = \{d^*, [H, d]\} = \Delta,
\]
because $[H, d]=d$ as Lefschetz theorem implies.

\hfill

\corollary\label{_product_harmo_Corollary_}
Let $(M, \omega)$ be a Riemannian manifold
equipped with a parallel differential form, and $\eta$ 
a harmonic form on $M$. Then $\omega\wedge \eta$ is
harmonic.

\hfill

{\bf Proof:} Follows from \ref{_kahler_ide_Proposition_}
(ii). \endproof

\hfill

This statement seems to be well known.

\hfill

Further on, we shall need the following trivial lemma.

\hfill

\lemma\label{_d_c_vanishes_on_harmonic_Lemma_}
Let $(M, \omega)$ be a compact Riemannian manifold
equipped with a parallel differential form, and $\eta$ 
a harmonic form on $M$. Consider the twisted de Rham
operator $d_c$ constructed above. Then $d_c(\eta)=0$.

\hfill

{\bf Proof:} Since $M$ is compact, $d^*\eta=0$. Then
$d_c\eta = d^*L_\omega\eta$. On the other hand,
$L_\omega\eta$ is harmonic, by
\ref{_product_harmo_Corollary_},
hence satisfies $d^*L_\omega\eta=0$.
\endproof

\subsection{The differential graded algebra
$(\ker d_c, d)$}

Let $(M,\omega)$ be a Riemannian manifold
equipped with a parallel form, and $d_c$ the twisted
de Rham operator constructed above. By construction,
$d_c$ is a differentiation of $\Lambda^*(M)$. Therefore,
$\ker d_c\subset \Lambda^*(M)$ is a subalgebra.
Since $d$ and $d_c$ supercommute, $d$ acts on
$\ker d_c$. We consider $(\ker d_c, d)$ as a differential
graded algebra (a DG-algebra).

\hfill

Recall that a homomorphism of DG-algebras is called
{\bf a quasi-\-iso\-mor\-phism} if it induces isomorphism on cohomology.

\hfill

\proposition\label{_ker_quasi-iso_Proposition_}
Let $(M,\omega)$ be a compact Riemannian manifold
equipped with a parallel form. Consider the natural embedding 
\begin{equation}\label{_ker_embedding_Equation_}
(\ker d_c, d)\hookrightarrow (\Lambda^*(M), d).
\end{equation}
Then this map is a quasi-isomorphism.

\hfill

{\bf Proof:} 
Let $\Lambda^*(M)_\alpha$ be the eigenspace of $\Delta$,
corresponding to the eigenvalue $\alpha$.  Since $\Delta$ is a self-adjoint
operator with discrete spectrum, we have a decomposition
$\Lambda^*(M) \cong \oplus_{\alpha}\Lambda^*(M)_\alpha$.
Consider the subcomplex
\begin{equation} \label{_d_on_eigenspa_Equation_}
   ... \stackrel d \arrow \Lambda^*(M)_\alpha \stackrel d \arrow
   \Lambda^{*+1}(M)_\alpha\stackrel d \arrow ...
\end{equation}
corresponding to an eigenvalue $\alpha$. Clearly,
for $\alpha\neq 0$, the complex
\eqref{_d_on_eigenspa_Equation_} is exact.
Let 
\begin{equation} \label{_d_on_eigenspa_ker_d_c_Equation_}
  ... \stackrel d \arrow (\ker d_c)_\alpha \stackrel d \arrow
   (\ker d_c)_\alpha\stackrel d \arrow ...
\end{equation}
be the action of $d$ on the $\alpha$-eigenspace of $\Delta$ on $(\ker d_c)$
($\Delta$ commutes with $d_c$ as 
\ref{_kahler_ide_Proposition_} implies). 

For $\alpha=0$, $(\ker d_c)_\alpha=\Lambda^*(M)_\alpha={\cal H}^*(M)$
as \ref{_d_c_vanishes_on_harmonic_Lemma_}
implies. To prove \ref{_ker_quasi-iso_Proposition_}
we need only to show that \eqref{_d_on_eigenspa_ker_d_c_Equation_}
has zero cohomology for $\alpha>0$. However, for any
closed form $\eta \in (\ker d_c)_\alpha$,
we have 
\[ \eta = \frac 1 \alpha (dd^* + d^* d)\eta =
    \frac 1 \alpha dd^*\eta
\]
and $d^*\eta$ lies inside $(\ker d_c)_\alpha$
as $d_c$ and $d^*$ commute (\ref{_kahler_ide_Proposition_}).
Therefore, $\eta$ is exact. This proves 
\ref{_ker_quasi-iso_Proposition_}.
\endproof

\hfill

The following claim is clear, as $\Delta_c$ and $\Delta$
commute, and 
$\{d_c, d_c^*\}^* = \{d_c^*, d_c\}=(-1)^{1-\tilde d_c}\{d_c, d_c^*\}.$

\hfill

\claim\label{_finite_Delta_c_Claim_}
Let $(M,\omega)$ be a compact Riemannian manifold
equipped with a parallel form, and 
$\Delta_c= \{d_c, d^*_c\}$  the operator constructed above.
Let $\Lambda^*(M)_\alpha$ be the eigenspace
of the Laplacian of eigenvalue $\alpha$.
Then $\Delta_c$ preserves $\Lambda^*(M)_\alpha$
and acts on $\Lambda^*(M)_\alpha$ as a self-adjoint
or anti-self-adjoint operator. In particular, $\Delta_c$
is diagonalizable, on some dense 
subspace of $\Lambda^*(M)\otimes_\R \C$

\endproof

\hfill

\remark
Notice that $\Delta_c$ is not a priori elliptic,
hence it has no spectral decomposition. However, it 
perserves the finite-dimensional eigenspaces of the
Laplacian, and is diagonalizable on these
eigenspaces.

\subsection{Pseudocohomology of the operator $d_c$}

\lemma \label{_im_d_c_ideal_Lemma_}
Let $(M,\omega)$ be a compact Riemannian manifold
equipped with a parallel form, and
$(\ker d_c, d)$ the differential
graded algebra constructed above. 
Consider the subspace
\begin{equation}\label{_ideal_d_c_Equation_}
V= (\ker d_c) \cap d_c(\Lambda^*(M)) \subset (\ker d_c).
\end{equation}
Then $V$ is a differential ideal in the
differential graded algebra $(\ker d_c, d)$.
In other words, $\ker d_c\cdot V \subset V$ and
$dV \subset V$.

\hfill

{\bf Proof:} 
Given $x \in \ker d_c$, $y\in V$,
$y = d_c z$, we write 
\[ d_c(x\wedge z) = (-1)^{\tilde d_c \tilde x}x\wedge d_c z.
\]
Therefore, $V$ is an ideal. To prove that $dV\subset V$,
we write $v\in V$ as $d_c(w)$, then 
$dv = (-1)^{\tilde d_c}d_c d w$.
\endproof

\hfill

\definition\label{_pseudoco_Definition_}
The quotient 
$\frac{(\ker d_c)}{(\ker d_c) \cap (\im d_c)}$
is called {\bf the pseu\-do-\-co\-ho\-mo\-lo\-gy} of $d_c$. As 
\ref{_im_d_c_ideal_Lemma_} implies, pseudo-cohomology
is a differential graded algebra. We denote it by
$(H^*_{c}(M), d)$. 

\hfill

\remark
We don't call $H^*_c(M)$ {\it cohomology} of $d_c$,
because $d_c^2$ is not necessarily zero. 
In the literature, the pseudo-cohomology 
of an operator is known under the name
{\bf twisted cohomology} (see e.g. in 
\cite{_Vaisman:twisted_}).

\hfill

\definition
Let $\eta\in \Lambda^*(M)$ be a form which satisfies 
$d_c\eta=d^*_c \eta=0$. Then $\eta$ is called {\bf
pseudo-harmonic}. The space of all pseudo-harmonic
forms is denoted by ${\cal H}^*_c(M)$. By 
\ref{_kahler_ide_Proposition_} (i), the de Rham
differential preserves ${\cal H}^*_c(M)$.

\hfill

\remark\label{_harmo_pseudoharmo_Remark_}
From \ref{_d_c_vanishes_on_harmonic_Lemma_}
it follows immediately that all harmonic forms
are pseudo-harmonic:  
${\cal H}^*(M)\subset {\cal H}^*_c(M)$.

\hfill

\proposition\label{_pseudoha_pseudoco_Proposition_}
Let $(M,\omega)$ be a compact Riemannian manifold
equipped with a parallel form, and 
\begin{equation}\label{_proje_to_pseu_Equation_}
  {\cal H}^*_c(M)\stackrel i \arrow H^*_{c}(M)
\end{equation}
the natural projection map. Then $i$ is an 
isomorphism, compatible with the de Rham differential.

\hfill

{\bf Proof:} 
We represent $\Lambda^*(M)$ as a (completion of)
a direct sum of eigenvalues of the Laplacian.
Using  $d_c$, $d_c^*$-invariance of these eigenspaces,
we may work with the associated decompositions within
these eigenspaces.  Abusing the language,
we approach $\Lambda^*(M)$ as if it 
were finite-dimensional, but in fact we work with
these eigenspaces, which are finite-dimensional.

From 
\[
(\Delta_c\alpha, \alpha)= (d_c\alpha, d_c\alpha) + (d_c^*\alpha, d_c^*\alpha)
\]
we obtain that $\ker\Delta_c= \ker d_c \cap \ker d_c^*$.
From $(d_c\alpha, \beta) = (\alpha, d_c^*\beta)$,
we find that $\ker d_c = (\im d_c^*)^\bot$,
$\ker d_c^* = (\im d_c)^\bot$, where $(\cdots)^\bot$
denotes the orthogonal complement. Therefore, 
\[ 
  \ker \Delta_c= (\im d_c)^\bot \cap (\im d_c^*)^\bot = 
(\im  d_c +\im d_c^*)^\bot.
\]
Given $\alpha\in \Lambda^*(M)$, let $\Pi\alpha$ denote the
orthogonal projection of $\alpha$ to $\ker \Delta_c$. Then
$\alpha-\Pi\alpha$ is orthogonal to $\ker \Delta_c$, hence
\begin{equation}\label{_proj_in_im_+im_Equation_}
  \alpha-\Pi\alpha\in  \bigg(\im  d_c +\im d_c^*\bigg). 
\end{equation}
Now assume that $\alpha\in \ker d_c$. The form
$\Pi(\alpha)$ also lies in $\ker d_c$,
because $\ker \Delta_c\subset \ker d_c$.
Therefore, $\alpha-\Pi\alpha$ lies in $\ker d_c$,
hence, is orthogonal to $\im d_c^*$.  Using
\eqref{_proj_in_im_+im_Equation_}, we obtain that
$\alpha-\Pi\alpha\in \im d_c$.

Therefore,
\begin{equation}\label{_ker_d_c_decompo_Equation_}
\ker d_c = (\ker d_c) \cap (\im d_c) \oplus {\cal H}^*_c(M).
\end{equation}
From \eqref{_ker_d_c_decompo_Equation_}
\ref{_pseudoha_pseudoco_Proposition_} 
follows directly.
\endproof

\hfill

\proposition\label{_H_c^*_quasi-iso_Proposition_}
Let $(M,\omega)$ be a compact Riemannian manifold
equipped with a parallel form, and
\begin{equation}\label{_ker_dc_to_pseudoco_Equation_}
(\ker d_c, d)\stackrel \pi \arrow (H^*_{c}(M), d)
\end{equation}
the homomorphism of differential
graded algebras constructed above.
Then $\pi$ is a quasi-isomorphism.

\hfill

{\bf Proof:} By definition, 
\eqref{_ker_dc_to_pseudoco_Equation_} is surjective.
To show that it is a quasi-isomorphism, we need to prove
that any $d$-closed $\eta \in\ker \pi$
is $d$-exact. However, $\ker\pi \subset (\ker d_c) \cap (\im d_c)$,
and by \eqref{_ker_d_c_decompo_Equation_} this space is
orthogonal to ${\cal H}^*_c(M)$.
Using \ref{_harmo_pseudoharmo_Remark_}
we obtain that any $\eta \in \ker\pi$ is orthogonal to
the space of harmonic forms. Using the spectral
decomposition, we obtain that
$\eta = \sum \eta_{\alpha_i}$, where
$\Delta\eta_{\alpha_i} = \alpha_i\eta_{\alpha_i}$,
and $\{\alpha_i\}$ are positive real numbers.
Since $\Delta$ commutes with $d_c$ and $d_c^*$,
the components $\eta_{\alpha_i}$
also belong to $\ker\pi$.
This gives $\eta_{\alpha_i}= \frac1 {\alpha_i} d d^*
\eta_{\alpha_i}$,
hence all the components $\eta_{\alpha_i}$
are $d$-exact. We obtain that $\eta$ is $d$-exact.
\ref{_H_c^*_quasi-iso_Proposition_} is proven.
\endproof

\hfill

\remark\label{_ddc_proof_Remark_}
The standard (and completely formal) agrument is used
to produce the $dd_c$-lemma from \ref{_H_c^*_quasi-iso_Proposition_}.
Let $\eta$ be a $d_c$-exact, $d$-, $d_c$-closed form on $M$.
We need to show that $\eta=dd_c\xi$. By definition,
$\eta$ represents 0 in $H^*_c(M)$.
Since $(\ker d_c, d)$ is quasi-isomorphic to $(H^*_c(M), d)$,
$\eta$ represents zero in the cohomology of 
$(\ker d_c, d)$. Therefore, $\eta = d \nu$, for
some $\nu\in \ker d_c$. Now, the class $[\nu]$ of 
$\nu$ in $H^*_c(M)$ satisfies $d[\nu]=0$.
Using \ref{_H_c^*_quasi-iso_Proposition_} again, we find that
$[\nu] - [\nu']= 0$, for some $d$-closed
form $\nu'\in \ker d_c$. Therefore, $\nu-\nu'= d_c\xi$.
Since $d\nu'=0$, this gives $d d_c \xi=d\nu=\eta$.

\hfill

\definition\label{_formal_DG_Definition_}
Let $(A^*, d)$, $(B^*, d)$ 
be graded commutative differential graded algebras
(DG-algebras, for short).
If $(A^*, d)$ and $(B^*, d)$ can be connected by a
sequence of quasi-isomorphisms 
\[ (A^*, d)\arrow (A^*_1, d),\ \ 
   (A_2^*, d_2)\arrow (A^*_1, d), \ \  ... \ \ (A_n^*, d_n)\arrow (B^*,d),
\]
the DG $(A^*, d)$ and $(B^*, d)$ are called {\bf weak 
equivalent}. A DG-algebra is called {\bf formal}
if it is weak equivalent to a DG-algebra with $d=0$.

\hfill

\corollary \label{_homotopi_equiv_for_pseudo_Corollary_}
Let $(M,\omega)$ be a compact Riemannian manifold
equipped with a parallel form, and $(H_c^*(M), d)$
its pseudohocomology DG-algebra. Then $(\Lambda^*(M), d)$
is weak equivalent to $(H_c^*(M), d)$. Moreover,
if every pseudoharmonic form is harmonic,
then $(\Lambda^*(M), d)$ is formal.

\hfill

{\bf Proof:} By \ref{_ker_quasi-iso_Proposition_},
the DG-algebra 
$(\Lambda^*(M), d)$ is quasi-\-iso\-mor\-phic to 
$(\ker d_c, d)$. By
\ref{_H_c^*_quasi-iso_Proposition_}, the DG-algebra
$(\ker d_c, d)$ is quasi-\-iso\-mor\-phic to $(H_c^*(M), d)$.
Finally, if all pseudoharmonic forms are harmonic,
the differential $d$ vanishes on ${\cal H}_c^*(M)$, 
and \ref{_pseudoha_pseudoco_Proposition_}
implies that $d=0$ on $(H_c^*(M), d)$.
\endproof

\hfill

\remark 
When $(M, \omega)$ is a compact K\"ahler manifold, 
$\Delta=\Delta_c$ as the K\"ahler identities imply.
In this situation, pseudoharmonic forms are
the same as harmonic. This implies the
celebrated result of \cite{_DGMS:Formality_}: 
for any compact K\"ahler manifold, its de Rham DG-algebra 
is formal.


\section{Structure operator for holonomy $G_2$-manifolds}


\subsection{$G_2$-manifolds}
\label{_G_2_Subsection_}

We base our exposition on \cite{_Hitchin:3-forms_}.

\hfill

\claim
Consider the natural action of $GL(7, \R)$ on the space
$\Lambda^3(V^*)$ of 3-forms on $V$, where $V= \R^7$. Then $GL(7, \R)$
acts on $\Lambda^3(V^*)$ with two open orbits.

{\bf Proof:} Well known (see e.g. \cite{_Salamon:holonomy_}).
\endproof

\hfill

\definition
A 3-form $\omega$ on $V= \R^7$.
is called {\bf non-degenerate} if it lies in an open
orbit. 

\hfill

The group $GL(7, \R)$ is 49-dimensional, and 
dimension of $\Lambda^3(V^*)$ is 35. Therefore,
a stabilizer of a non-degenerate 3-form has 
dimension 14. This stabilizer is a Lie group, of 
dimension 14, called $G_2$. For one orbit it is a compact
form of $G_2$, for another orbit a non-compact 
real form. We call a non-degenerate 3-form 
$\omega$ on $V= \R^7$ {\bf positive} if its
stabilizer is a compact form of $G_2$.

\hfill

Given a 3-form
$\omega\in\Lambda^3(V^*)$, consider an
$\Lambda^7(V^*)$-valued scalar product
$V\times V\arrow \Lambda^7(V^*)$,
\[ x, y \stackrel {\tilde g}\arrow \frac 1 6
   (\omega\cntrct x)\wedge (\omega\cntrct y)\wedge\omega.
\]
It is easy to check that ${\tilde g}$ is non-degenerate
when $\omega$ is non-degenerate, and sign-definite
when $\omega$ is positive. Consider $\tilde g$
as a section of $V^*\otimes V^* \otimes \Lambda^7(V^*)$,
and denote by $K$ its determinant, $K\in \Lambda^7(V^*)^9$.
Since 9 is odd, $K$ gives an orientation on $V$.
Let $k:= \sqrt[9]{K}$ be the corresponding
section of $\Lambda^7(V^*)$, and $g:= k^{-1}\tilde g$
the $\R$-valued bilinear symmetric form associated 
with $\tilde g$. Assume that $\omega$ is positive.
A direct calclulation implies
that $g$ is positive definite, and in some
orthonormal basis $e_1, ..., e_7\in V^*$, $\omega$
is written as
\begin{equation}\label{_omega_in_basis_Equation_}
\begin{aligned}
\omega &= (e_1\wedge e_2 + e_3\wedge e_4)\wedge e_5 + 
 (e_1\wedge e_3 - e_2\wedge e_4)\wedge e_6 \\ + &
 (e_1\wedge e_4 - e_2\wedge e_3)\wedge e_7 + 
 e_5\wedge e_6\wedge e_7.
\end{aligned}
\end{equation}

\hfill

\definition
Let $M$ be a 7-dimensional smooth manifold, and
$\omega\in \Lambda^3(M)$ a 3-form. $(M, \omega)$ is called
{\bf a $G_2$-manifold} if $\omega$ is non-degenerate
and positive everywhere on $M$. We consider $M$
as a Riemannian manifold, with the Riemannian structure
determined by $\omega$ as above. The manifold 
$(M, g,\omega)$ is called {\bf a holonomy $G_2$-manifold}
if $\omega$ is parallel with respect to the Levi-Civita
connection associated with $g$. Further on, we
shall consider only holonomy $G_2$ manifolds, and
(abusing the language) omit the word ``holonomy''.

\hfill

\remark
Holonomy $G_2$-manifolds have 
long and distinguished history.
They appear in M. Berger's list of
irreducible holonomies (\cite{_Berger:list_}). 
Local examples of holonomy $G_2$-manifolds were unknown 
untill R. Bryant's work of mid-1980-ies
(\cite{_Bryant:holonomy_constru_}).
Then R. Bryant and S. Salamon constructed
a complete examples of holonomy $G_2$-manifoldy
(\cite{_Bryant_Salamon_}), and D. Joyce
(\cite{_Joyce_G2_}) constructed and studied 
compact holonomy $G_2$-manifolds at great 
length. For details of D. Joyce's 
construction, see \cite{_Joyce_Book_}.
Since then, the $G_2$-manifolds become
crucially important in many areas
of string physics, especially 
in M-theory.

\hfill

Under the $G_2$-action, the space
$\Lambda^*(M)$ splits into irreducible
representations, as follows.

\begin{equation}\label{_Lambda^*_splits_Equation_}
\begin{aligned}
\Lambda^2(M)&
\cong \Lambda^2_7(M)\oplus \Lambda^2_{14}(M), \\
\Lambda^3(M) &\cong
\Lambda^3_1(M)\oplus \Lambda^3_7(M)\oplus \Lambda^3_{27}(M)
\end{aligned}
\end{equation}
where $\Lambda^i_j(M)$ is an irreducible $G_2$-representation
of dimension $j$. Clearly, $\Lambda^*(M)\cong
\Lambda^{7-*}(M)$ as a $G_2$-representation,
and the spaces $\Lambda^4(M)$, $\Lambda^5(M)$
split in a similar fashion. The spaces
$\Lambda^0$, $\Lambda^1$ are irreducible.

The spaces $\Lambda^i_j(M)$ are defined explicitly, 
in a following way. $\Lambda^2_7(M)$ is
$\Lambda_{*\omega}(\Lambda^6(M))$, where 
$\Lambda_{*\omega}$ is the Hermitian adjoint
to $L_{*\omega}(\eta) = *\omega\wedge \eta$
(see Section \ref{_para_form_Section_}).  
The space $\Lambda^2_{14}(M)$ is identified with 
${\goth g}_2\subset \goth{so}(TM)$ under the 
standard identification $\Lambda^2(M)=\goth{so}(TM)$.
The space $\Lambda^3_1(M)$ is generated by $\omega$,
$\Lambda^3_7(M)$ is equal to
$\Lambda_\omega(\Lambda^6(M))$, where
$\Lambda_\omega$ is the Hermitian adjoint
of $L_\omega(\eta) = \omega\wedge \eta$
(see Section \ref{_para_form_Section_}).
Finally, $\Lambda^3_{27}(M)$ is identified with
$(\ker L_\omega)\cap (\ker \Lambda_\omega)\subset \Lambda^3(M)$.

\hfill

\remark
Notice that the operators $C$, $L_\omega$, 
$\Lambda_\omega$ from Section \ref{_para_form_Section_}
are clearly $G_2$-invariant.

\hfill

From the construction, it is clear that the splitting
\eqref{_Lambda^*_splits_Equation_} can be obtained 
via the operators $L_\omega$, $\Lambda_\omega$, 
$L_{*\omega}$, $\Lambda_{*\omega}$.
By \ref{_kahler_ide_Proposition_}
these operators commute with the Laplacian.
Therefore, harmonic forms also split:
 
\begin{equation}\label{_harmo_splits_Equation_}
\begin{aligned}
{\cal H}^2(M)&
\cong {\cal H}^2_7(M)\oplus {\cal H}^2_{14}(M), \\
{\cal H}^3(M) &\cong
{\cal H}^3_1(M)\oplus {\cal H}^3_7(M)\oplus {\cal H}^3_{27}(M)
\end{aligned}
\end{equation}
and similar splitting occurs on ${\cal H}^4(M)$ and 
${\cal H}^5(M)$.

The following result is well known
and is implied by a Bochner-Lich\-ne\-ro\-wicz-\-type argument using 
Ricci-flatness of holonomy $G_2$-manifolds.

\hfill

\claim\label{_H_7_parallel_Claim_}
Let $M$ be a compact $G_2$-manifold,
and $\eta\in {\cal H}^i_7(M)$ a harmonic form. Then
$\eta$ is parallel. Moreover, if $H^1(M)=0$, then 
${\cal H}^i_7(M)=0$ ($i=1, 2, 3, 4, 5, 6$).

{\bf Proof:} See \cite{_Joyce_Book_}.
\endproof

\hfill

\remark
A $G_2$-manifold is Ricci-flat,
as shown by E. Bonan (\cite{_Bonan:G_2_}).
Then $\pi_1(M)$ is finite, unless $M$ has a finite
covering which is isometric to $T\times M'$,
where $M'$ is a manifold with special holonomy,
and $T$ a torus. When $\pi_1(M)$ is finite,
${\cal H}^i_7(M)=0$ as \ref{_H_7_parallel_Claim_}
implies. 

\hfill

We shall also need the following linear-algebraic result,
which is well known.
Let $M$ be a $G_2$-manifold, and 
\[ \Lambda^2(M) \cong \Lambda^2_7(M)\oplus \Lambda^2_{14}(M)
\]
the decomposition defined above. Consider the
operator \[ * \circ L_\omega:\; \Lambda^2(M)\arrow \Lambda^2(M).\]
This operator is $G_2$-invariant, hence by Schur's lemma
acts on $\Lambda^2_7(M)$ and $\Lambda^2_{14}(M)$
as scalars. These scalars are computed as follows

\hfill

\claim\label{_*_on_Lambda^2_Claim_}
For any  $\alpha \in \Lambda^2_7(M)$, we have
$*  L_\omega\alpha = 2 \alpha$. For
$\alpha \in \Lambda^2_{14}(M)$, we have
$*  L_\omega\alpha = -\alpha$.

\hfill

{\bf Proof:} See e.g. \cite{_Bryant:G2_}, (2.32). \endproof

\subsection{Structure operator for  $G_2$-manifolds}

Let $(M, \omega)$ be a $G_2$-manifold. 
We have two parallel forms on $M$: $\omega$ and
$*\omega$, and the results of Section 
\ref{_para_form_Section_} can be applied to
$\omega$ and $*\omega$ as well. 

We denote by $C$, $C_{*\omega}$ the corresponding
structure operators, and by $d_c$ the operator
$\{ C, d\}$. 

This part of the paper is a
pure linear algebra. We never use 
the holonomy property: throughout this subsection,
there is no need to assume that our $G_2$-manifold
has holonomy in $G_2$.

Consider the operator $C^2= \frac 1 2 \{ C,C\}$.
Being a supercommutator of two differentiations,
this operator is a differentiation.

\hfill

\claim\label{_C^2_Claim_}
Under these assumptions,
\begin{equation}\label{_C^2_Equation_}
C^2 = 3 C_{*\omega}
\end{equation}

{\bf Proof:}
Both sides of \eqref{_C^2_Equation_}
are differentiations, and vanish on $\Lambda^0(M)$.
Therefore, to prove \eqref{_C^2_Equation_} it suffices
to check that $C^2 = 3 C_{*\omega}$ on $\Lambda^1(M)$.
Both $C^2$ and $C_{*\omega}$ define $G_2$-invariant
map from $\Lambda^1(M)$ to $\Lambda^3_7(M)$.
By Schur's lemma, these operators are proportional.
To show that the coefficient of proportionality is 3,
we compute $C^2$ and $C_{*\omega}$ on $e_1$, using
\eqref{_omega_in_basis_Equation_}. \endproof

\hfill

A similar argument gives the following claim

\hfill

\claim
Under the above assumptions, we have
\begin{equation}\label{_L_omega_C^*_commu_Equation_}
\{L_\omega, C^*\} = -3 C_{*\omega}.
\end{equation}

{\bf Proof:} The operator $C^*$ takes a form
\begin{equation}\label{_C^*_coordinates_Equation_}
\begin{aligned}
\  & C^*(e_{i_1}\wedge e_{i_2}\wedge ...) \\ &= 
\sum_{k_1<k_2} (-1)^{(i_{k_1}-1)i_{k_2}} 
C^*(e_{i_{k_1}}\wedge e_{i_{k_2}})
\wedge e_{i_1}\wedge e_{i_2}\wedge ...\wedge \check e_{i_{k_1}}\wedge ... \wedge \check e_{i_{k_2}}\wedge ...
\end{aligned}
\end{equation}
where $C^*(e_{i_{k_1}}\wedge e_{i_{k_2}})$
is the usual crossed product of vectors $e_{i_{k_1}}, e_{i_{k_2}}$
on the space equipped with a 3-form and a non-degenerate
bilinear symmetric form. From \eqref{_C^*_coordinates_Equation_}
it is clear that $C^*$ is a second order differential
operator on the algebra $\Lambda^*(M)$ (differential
operators on a graded commutative algebra are understood
in the sense of Grothendieck - see e.g. \cite{_Verbitsky:HKT_}).
Then $\{L_\omega, C^*\}$ is a first order differential
operator. An elementary calculation gives $C^*\omega=0$.
Therefore, $\{L_\omega, C^*\}$ is a differentiation.
To compare $\{L_\omega, C^*\}$ with $-3 C_{*\omega}$,
we need to check that $\{L_\omega, C^*\} = -3 C_{*\omega}$
on $\Lambda^1(M)$. Both of these operators are $G_2$-invariant,
and Schur's lemma implies that they are proportional
on $\Lambda^1(M)$. To compute the coefficient of
proportionality, it suffices to compute 
$\{L_\omega, C^*\}$, $C_{*\omega}\}$
on some vector, e.g. $e_1$.
\endproof

\hfill

\claim\label{_C_on_Lambda^3_Claim_}
Under the above assumptions, $C:\; \Lambda^3(M)\arrow \Lambda^4(M)$
is an isomorphism. Moreover, $C\omega = 2 *\omega$.

\hfill

{\bf Proof:} 
Clearly, $C$ preserves the decomposition of $\Lambda^*(M)$
onto $G_2$-invariant summands as in 
\eqref{_Lambda^*_splits_Equation_}.
We write $\omega$ in orthonormal basis as in 
\eqref{_omega_in_basis_Equation_}. The equation 
$C\omega = 2 *\omega$ is given by a direct calculation.
Given a 3-form $\theta \in \Lambda^3_7(M)$
and applying \eqref{_L_omega_C^*_commu_Equation_},
we obtain $\Lambda^*(C\theta)= -3 (C^*_{*\omega})\theta$.
However, $C^*_{*\omega}:\; \Lambda^3_7(M)\arrow
\Lambda^1(M)$ is an isomorphism, because
$C_{*\omega}:\; \Lambda^1(M)\arrow \Lambda^3_7(M)$
is non-zero. To prove \ref{_C_on_Lambda^3_Claim_},
it remains to show that $C$ is an isomorphism on
$\Lambda^3_{27}(M)$. By Schur's lemma,
for this it suffices to show that
$C\restrict{\Lambda^3_{27}(M)}$ is non-zero.

Consider the form 
$\eta=e_5\wedge (e_1 \wedge e_2 - e_3 \wedge e_4)$.
Clearly, $\Lambda_\omega\eta=0$ and
$L_\omega\eta=0$. Therefore, $\eta \in \Lambda^3_{27}(M)$.

From \eqref{_C^*_coordinates_Equation_} we find
that $C^*(e_1 \wedge e_2 - e_3 \wedge e_4)=0$,
hence $e_1 \wedge e_2 - e_3 \wedge e_4$ lies in $\Lambda^2_{14}(M)$.
This gives
\begin{equation}\label{_C_eta_long_expli_Equation_}
\begin{aligned}
 C(\eta)=&C(e_5\wedge (e_1 \wedge e_2 - e_3 \wedge e_4))\\
 &=C(e_5) \wedge (e_1 \wedge e_2 - e_3 \wedge e_4)\\ &=
 (e_1 \wedge e_2 + e_3 \wedge e_4+ e_6 \wedge e_7)
 \wedge (e_1 \wedge e_2 - e_3 \wedge e_4) \\& = 
 e_6 \wedge e_7\wedge (e_1 \wedge e_2 - e_3 \wedge e_4)
\end{aligned}
\end{equation}
We obtain that $C(\eta)\neq 0$. 
\ref{_C_on_Lambda^3_Claim_} is proven.
\endproof

\hfill

\remark\label{_C_eta_Lambda^3_27_Remark_}
The calculation \eqref{_C_eta_long_expli_Equation_} gives
\begin{equation}\label{_C_on_Lambda^3_27_*_Equation_}
 C(\eta)=-*\eta 
\end{equation}
and by Schur's lemma this 
equation holds for all $\eta\in \Lambda^3_{27}(M)$.

\hfill

\proposition\label{_C_isomo_on_Lambda_7_Proposition_}
Let $(M,\omega)$ be a $G_2$-manifold, and $C$ its structure
operator. Then $C$ induces isomorphisms
\begin{equation}\label{_C_on_Lambda_7_Equation_}
\Lambda^i_7(M) \stackrel C \arrow \Lambda^{i+1}_7(M),
\end{equation}
$(i=1, 2, 3, 4, 5)$.

\hfill

{\bf Proof:}
By Schur's lemma, \eqref{_C_on_Lambda_7_Equation_}
is either an isomorphism or zero. For $i=1$,
$i=2$ \eqref{_C_on_Lambda_7_Equation_} is non-zero
as follows from \ref{_C^2_Claim_}. For $i=3$,
\eqref{_C_on_Lambda_7_Equation_} is non-zero by 
\ref{_C_on_Lambda^3_Claim_}. Using 
\[ 
  C(\phi \wedge \psi) = C(\phi) \wedge \psi + 
  (-1)^{\tilde \phi}\phi \wedge C(\psi),
\]
we find that $*C*$ is Hermitian adjoint to $C$. 
On the other hand, \eqref{_C_on_Lambda_7_Equation_}
is an isomorphism if and only if
\[ 
\Lambda^{i+1}_7(M) \stackrel {C^*} \arrow \Lambda^{i}_7(M),
\]
is an isomorphism. Using $C^*=*C*$,
we obtain that \ref{_C_isomo_on_Lambda_7_Proposition_}
$i=k$ is implied by \ref{_C_isomo_on_Lambda_7_Proposition_}
for $i=6-k$. Therefore, the already proven
assertions of \ref{_C_isomo_on_Lambda_7_Proposition_}
for $i=1,2,3$ imply \ref{_C_isomo_on_Lambda_7_Proposition_}
for $i=4,5$. \endproof


\section{Pseudocohomology for $G_2$-manifolds}


\subsection{De Rham differential on $\Lambda^*_7(M)$}

To study the pseudocohomology, we use the following well
known lemma (appearing in a different form in 
\cite{_Trieste_G2:Fernandez_Ugarte_}
and \cite{_Dolbeault_G2:Fernandez_Ugarte_}).

\hfill

\lemma\label{_d_on_Lambda_7_Lemma_}
Let $\eta \in \Lambda^k_7(M)$ be a differential form
on a holonomy $G_2$-manifold (not necessarily compact),
where $0<k<5$ is an integer. Fix parallel 
$G_2$-invariant isomorphisms
\begin{equation}\label{_iso_Lambda_7_Equation_}
\Lambda^k_7(M)\stackrel{\tau_{i,k}} \arrow \Lambda^i_7(M), 
\end{equation}
for all 
$i = 1, 2, 3, 4,5$ (by Schur's lemma, 
these isomorphisms are well defined, 
up to a constant). \footnote{Using 
\ref{_C_isomo_on_Lambda_7_Proposition_},
we could use the powers of $C$ to define
the isomorphisms \eqref{_iso_Lambda_7_Equation_}.}
Denote by 
$d_7:\; \Lambda^i_7(M)\arrow \Lambda^{i+1}_7(M)$
the $\Lambda^*_7$-part of the de Rham differential.
Then $d_7(\eta)=0$ if and only if $d_7(\tau_{k,i}\eta)=0$
for any $i=1,2,3, 4$.

\hfill

{\bf Proof:}
Consider the Levi-Civita connection 
\begin{equation}\label{_nabla_on_Lambda_7_Equation_}
\nabla:\; \Lambda^i_7(M)\arrow \Lambda^i_7(M)\otimes \Lambda^1_7(M).
\end{equation}
The operator $d_7$ is obtained as a composition of 
\eqref{_nabla_on_Lambda_7_Equation_} and a $G_2$-invariant
pairing $\Lambda^i_7(M)\otimes \Lambda^1_7(M)\arrow \Lambda^{i+1}_7(M)$.
Using an irreducible decomposition of 
$\Lambda^1_7(M)\otimes \Lambda^{i}_7(M)$
(see e.g. \cite{_Bryant:G2_}), we find that
$\Lambda^1_7(M)\otimes \Lambda^{i}_7(M)$
contains a unique irreducible summand isomorphic
to $\Lambda^*_7(M)$ as a $G_2$-representation.
It is clear that $d_7:\; \Lambda^i_7(M)\arrow \Lambda^{i+1}_7(M)$
is obtained as a composition of 
\eqref{_nabla_on_Lambda_7_Equation_} and the projection
to this $\Lambda^*_7(M)$-summand. Therefore, the
following diagram is commutative, up to a constant multiplier
\begin{equation}
\begin{CD}
\Lambda^k_7(M)@>{d_7}>> \Lambda^{k+1}_7(M)\\
@V{\tau_{i,k}}VV @V{\tau_{i+1,k+1}}VV\\
\Lambda^i_7(M)@>{d_7}>> \Lambda^{i+1}_7(M).
\end{CD}
\end{equation}
We obtain that $\tau_{i+1,k+1}d_7(\eta)=0$
if and only if $d_7(\tau_{i,k}\eta)=0$.
This proves \ref{_d_on_Lambda_7_Lemma_}.
\endproof

\subsection{Computations of pseudocohomology}

\theorem\label{_pseudoha_G_2_Theorem_}
Let $(M,\omega)$ be a compact $G_2$-manifold,
${\cal H}^*(M)$ the space of harmonic forms,
and ${\cal H}^*_c(M)\supset {\cal H}^*(M)$ 
the space of pseudoharmonic forms. Then

\begin{description}
\item[(i)] 
${\cal H}^i_c(M)={\cal H}^i(M)$ for
all $i\neq 3,4$. 

\item[(ii)] The
orthogonal complement\footnote{This
notation has the following meaning: ${\cal H}^i_c(M)_{>0}$
is a sum of all positive eigenspaces of Laplacian
acting on ${\cal H}^i_c(M)$.}
 ${\cal H}^i_c(M)_{>0}$
to ${\cal H}^i(M)$ in ${\cal H}^i_c(M)$ 
lies in $\Lambda^i_{27}(M)$.
\item[(iii)]
$*({\cal H}_c^3(M)_{>0})={\cal H}_c^4(M)_{>0}$.
Moreover, ${\cal H}_c^3(M)_{>0}$ is generated by
all solutions of the equation $d\eta= \mu *\eta$,
for all $\mu\in \C$, $\mu\neq 0$, $\eta \in \Lambda^3_{27}(M)$.
\end{description}

{\bf Proof:} Consider the orthogonal decomposition
${\cal H}^*_c(M)= {\cal H}^*(M)\oplus{\cal H}^*(M)_{>0}$.
Since $\Delta$ preserves ${\cal H}^*_c(M)$,
$\Delta$ acts diagonally on ${\cal H}^*_c(M)$,
and ${\cal H}^*(M)_{>0}$ is generated by
all eigenvectors of $\Delta$ with 
non-zero eigenvalue. Therefore,
$d$ preserves ${\cal H}^*(M)_{>0}$

By \ref{_homotopi_equiv_for_pseudo_Corollary_}
${\cal H}^*_c(M)$ is quasi-isomorphic to 
${\cal H}^*(M)$. Therefore, cohomology of $d$ on
${\cal H}^*(M)_{>0}$ is zero. Now,
\ref{_pseudoha_G_2_Theorem_} (i)
is implied by the following claim

\hfill

\claim\label{_exa_H^*_c_Claim_}
Let $(M,\omega)$ be a compact $G_2$-manifold,
and $\eta\in {\cal H}^i_c(M)$ a non-zero exact 
pseudo-harmonic form. Then $i=4$.

\hfill

{\bf Proof:}
To prove  \ref{_pseudoha_G_2_Theorem_} (i)
suffices to prove \ref{_exa_H^*_c_Claim_}
for $i\leq 4$. Indeed, this will imply
that ${\cal H}^i(M)_{>0}=0$ for $i<3$,
but the Hodge $*$-operator preserves 
${\cal H}^*_c(M)$, and exchanges
${\cal H}^i(M)_{>0}$ and ${\cal H}^{7-i}(M)_{>0}$,
hence ${\cal H}^i(M)_{>0}=0$ for $i=1,2$
implies ${\cal H}^i(M)_{>0}=0$ for $i=5,6$.

Now, \ref{_pseudoha_G_2_Theorem_} (i)
is equivalent to \ref{_exa_H^*_c_Claim_}
as we have shown above. The same argument shows that 
\ref{_exa_H^*_c_Claim_} for $i\leq 4$.
implies \ref{_pseudoha_G_2_Theorem_} (i)
and the full statement of \ref{_exa_H^*_c_Claim_}.

Let $\eta = d\alpha$ be a $d$-exact 1-form in 
${\cal H}^1_c(M)$, $\alpha \in {\cal H}^2_c(M)$. 
Then $C\eta = Cd\alpha = - dC\alpha =0$
(the middle equation is implied by $d_c\alpha=0$).
Therefore, $C\eta=0$. However, $C$ is clearly injective
on $\Lambda^1(M)$. This proves \ref{_exa_H^*_c_Claim_} for $i=1$.

Let now $\eta = d\alpha$ be a $d$-exact 2-form in 
${\cal H}^2_c(M)$, $\alpha \in {\cal H}^1_c(M)$.
Using $d_c\eta=0$, we obtain
\begin{equation}\label{_C_eta_exact_Equation_}
0 = \{d, c\} \alpha = C \eta + d C\alpha.
\end{equation}
Write the decomposition $\eta= \eta_7+\eta_{14}$
induced by 
$\Lambda^2(M) \cong \Lambda^2_7(M)\oplus \Lambda^2_{14}(M)$.
Then \eqref{_C_eta_exact_Equation_} gives
$dC\eta = dC\eta_7=0$. From \ref{_d_on_Lambda_7_Lemma_}
we infer that $d_7\eta_7=0$. Consider the top degree forms
\begin{equation}\label{_equa_eta_7_product_d_alpha_Equation_}
\eta_7\wedge d\alpha \wedge \omega = \eta_7\wedge \eta_7\wedge \omega.
\end{equation}
(the equality holds by Schur's lemma as $\eta_7$ is
the $\Lambda^2_7(M)$-part of $\eta=d\alpha$).
Since $\Lambda^2_7(M)$ is an irreducible representation
of $G_2$, by Schur's lemma the 2-form
$\eta_7 \arrow \int \eta_7\wedge \eta_7\wedge \omega$
is sign-definite (negative definite, as
\ref{_*_on_Lambda^2_Claim_} implies).
Then  $\int \eta_7\wedge \eta_7\wedge \omega<0$
unless $\eta_7=0$. However, by 
\eqref{_equa_eta_7_product_d_alpha_Equation_}
\[ \int \eta_7\wedge \eta_7\wedge \omega
= \int \eta_7\wedge d\alpha \wedge \omega 
= - \int d\eta_7\wedge \alpha \wedge \omega =
   \int d_7\eta_7\wedge \alpha \wedge \omega =0
\]
as $d_7\eta_7=0$. We obtain that
$\eta\in \Lambda^2_{14}(M)$. 
Using \ref{_*_on_Lambda^2_Claim_} again, we
obtain  that $\int \eta\wedge\eta\wedge \omega>0$
unless $\eta=0$. However, $\eta$ is exact, hence this
integral vanishes, bringing $\eta=0$. 
We proved \ref{_exa_H^*_c_Claim_} for $i=2$.

Now, let $\eta = d\alpha$ be a $d$-exact 3-form in 
${\cal H}^3_c(M)$, $\alpha \in {\cal H}^2_c(M)$.
To finish the proof of 
\ref{_exa_H^*_c_Claim_}, we need to show that
$\eta=0$. 

Since $d^*$ commutes with $d_c, d^*_c$,
we have $d^*\alpha\in {\cal H}^1_c(M)$. As we have shown
above,  ${\cal H}^1_c(M)= {\cal H}^1(M)$, and therefore
$d^*\alpha$ is harmonic. A $d^*$-exact harmonic form vanishes.
Therefore, $d^*\alpha=0$. 

Then $0 = d_c \alpha = d^* L_\omega \alpha$. Similarly,
$0 = d^*_c\alpha = \Lambda d\alpha$. 
Using 
\[ \Lambda^2(M) \cong \Lambda^2_7(M)\oplus \Lambda^2_{14}(M),
\]
write the decomposition
$\alpha = \alpha_7 + \alpha_{14}$.
Then $L_\omega \alpha = 2*\alpha_7 - *\alpha_{14}$
as follows from \ref{_*_on_Lambda^2_Claim_}.
Therefore, $d^*L_\omega \alpha = * d(2\alpha_7 - \alpha_{14})$.
We obtain that $\alpha$ satisfies the following:
\begin{equation}\label{_d_alpha_in_Lambda^2_Equation_}
d^*\alpha =0, \ \  d(2\alpha_7 - \alpha_{14})=0, \ \ \Lambda d\alpha=0.
\end{equation}
Clearly, $d^*\alpha=0$ is equivalent to $d*\alpha=0$.
Also, $\{L_\omega, d\}=0$ ($\omega$ is closed). 
Using $* \alpha_7 = \frac 1 2 \alpha_7\wedge \omega$, 
$* \alpha_{14} = -\alpha_{14}\wedge \omega$
(\ref{_*_on_Lambda^2_Claim_}), we rewrite $d*\alpha=0$ as 
$L_\omega(d \alpha_{14}-\frac 1 2 d\alpha_7)=0$.
From \eqref{_d_alpha_in_Lambda^2_Equation_} we obtain
$L_\omega(d \alpha_{14}- 2 d\alpha_7)=0$.
Comparing these equations, we find
\begin{equation}\label{_L_omega_d_alpha_7_14_Equation_}
L_\omega(d \alpha_{14})=0 , \ \ L_\omega(d \alpha_{7})=0
\end{equation}
Using \ref{_*_on_Lambda^2_Claim_} again, we find that
\eqref{_L_omega_d_alpha_7_14_Equation_} implies
$d^*\alpha_{14}=d ^*\alpha_{7}=0$.

Now, $C^*\Lambda^2_{14}(M)=0$ because $C^*$ is $G_2$-invariant.
Using $d^*\alpha= d^*\alpha_7=0$, we obtain
\[ 0 = d_c^*\alpha = \{d^*, C^*\}\alpha = 
   d^* C^* (\alpha_{14} +\alpha_7) = d^* C^* \alpha_7= d_c^*\alpha_7.
\]
This implies 
\begin{equation}\label{_d^*_c_alpha_i_Equation_}
d_c^*\alpha_7=d_c^*\alpha_{14}=0
\end{equation}
Applying $d_c^*=\{d, \Lambda_\omega\}$, 
we find that \eqref{_d^*_c_alpha_i_Equation_} brings
\begin{equation}\label{_d^*_Lambda_alpha_i_Equation_}
\Lambda_\omega d\alpha_7=\Lambda_\omega d\alpha_{14}=0.
\end{equation}
Comparing \eqref{_d^*_Lambda_alpha_i_Equation_} and
\eqref{_L_omega_d_alpha_7_14_Equation_}, we find that
\begin{equation}\label{_d_alpha_in_Lambda_27_Equation_}
d\alpha_7, d\alpha_{14}\in \Lambda^3_{27}(M).
\end{equation}

This gives $\eta=d\alpha\in \Lambda^3_{27}(M)$.
Since $d_c\eta=0$, we have $dCd\alpha=0$, and the form 
$C\eta=Cd\alpha$ is closed. Therefore, 
\begin{equation}\label{_eta_C_eta_int_Equation_}
\int \eta \wedge C\eta = \int d\alpha\wedge C d\alpha =0.
\end{equation}
However, on $\Lambda^3_{27}(M)$, the form
$\eta \arrow \int \eta \wedge C\eta$ is non-zero
(\ref{_C_on_Lambda^3_Claim_}), hence, by Schur's lemma, 
sign-definite.\footnote{From 
\ref{_C_eta_Lambda^3_27_Remark_} it follows
that this form is negative definite.}
Therefore, \eqref{_eta_C_eta_int_Equation_}
implies that $\eta=0$. This proves \ref{_exa_H^*_c_Claim_}
for $i=3$. We finished the proof of \ref{_exa_H^*_c_Claim_}.
The proof of \ref{_pseudoha_G_2_Theorem_} (i) is also finished.
\endproof

\hfill

Let $\alpha \in \Lambda^3(M)$, and 
$\alpha = \alpha_1 + \alpha_7 +\alpha_{27}$
its decomposition induced by \eqref{_Lambda^*_splits_Equation_}.
To prove \ref{_pseudoha_G_2_Theorem_} (ii), we use the following
trivial observation:
\begin{equation}\label{_alpha_1_7_via_Lambda_Equation_}
\alpha_1 = \frac 1 7 L_\omega \Lambda_\omega \alpha, \ \ 
\alpha_7 = \frac 1 4  \Lambda_\omega L_\omega\alpha.
\end{equation}
Similarly, for $\eta \in \Lambda^4(M)$, 
$\eta = \eta_1 + \eta_7 +\eta_{27}$, we have
\begin{equation}\label{_eta_1_7_via_Lambda_Equation_}
\eta_1 = \frac 1 7  \Lambda_\omega L_\omega \eta, \ \ 
\eta_7 = \frac 1 4  L_\omega\Lambda_\omega \eta.
\end{equation}

Assume now that $\alpha \in {\cal H}^3_c(M)_{>0}.$
Then $d^*\alpha=0$ as \ref{_pseudoha_G_2_Theorem_} (i)
implies. Therefore
\[
0= d_c \alpha =\{ L_\omega , d^*\} \alpha = d^* L_\omega \alpha.
\]
From \eqref{_alpha_1_7_via_Lambda_Equation_},
we obtain
\begin{equation}\label{_d^*_alpha_7_Equation_}
d^*\alpha_7 = \frac 1 4 d^* \Lambda_\omega L_\omega \alpha=
-\Lambda_\omega d_c \alpha=0
\end{equation}
This implies
\begin{equation}\label{_d^*_c_alpha_7_Equation_}
d^*_c\alpha_7 = \{  d^*, C^*\} \alpha_7 = d^* C^* \alpha
\end{equation}
(the last equation holds because 
\[ \ker C^*\restrict{\Lambda^3(M)}= \Lambda^3_1(M)\oplus \Lambda^3_{27}(M)
\]
as $G_2$-decomposition implies).
However,
\[
d^* C^* \alpha = d^*_c \alpha =0
\]
since $d^*\alpha=0$. Then \eqref{_d^*_c_alpha_7_Equation_}
gives $d^*_c\alpha_7=0$. Similarly, 
\begin{equation}\label{_d^_c_alpha_7_Equation_}
d_c\alpha_7 = d^* L_\omega \alpha_7 = d^* L_\omega \alpha
\end{equation}
(here we use 
\[ \ker L_\omega\restrict{\Lambda^3(M)}= \Lambda^3_1(M)\oplus \Lambda^3_{27}(M)
\]
also implied by $G_2$-decomposition).
Using 
\[ 0 = d_c \alpha = \{ d^*, L_\omega\}\alpha = d^*L_\omega \alpha,
\]
we infer from \eqref{_d^_c_alpha_7_Equation_}
$d_c\alpha_7=0$. This gives 
$\alpha_7 \in {\cal H}^3_c(M)_{>0}.$

Now, by \eqref{_alpha_1_7_via_Lambda_Equation_},
\begin{equation}\label{_d_L_Lambda_alpha_7_Equation_}
0= d L_\omega \Lambda_\omega \alpha_7 = 
L_\omega \Lambda_\omega d\alpha_7 + L_\omega d^*_c\alpha_7=
L_\omega \Lambda_\omega d\alpha_7.
\end{equation}
Using \eqref{_eta_1_7_via_Lambda_Equation_},
we obtain that \eqref{_d_L_Lambda_alpha_7_Equation_}
gives $d\alpha_7 \in \Lambda^4_1(M)$. 
This means that $d\alpha_7 = f* \omega$, 
where $f\in C^\infty(M)$ is a function. Therefore,
$0 = d^2 \alpha_7 = df\wedge *\omega$. This leads
to $df=0$, as the map 
\[ \Lambda^1(M) \stackrel{L_{*\omega}}\arrow \Lambda^5(M) 
\]
is clearly injective. Therefore, $\alpha_7$ is harmonic,
hence $\alpha_7=0$. 

We have shown that 
\[ 
{\cal H}^3_c(M)_{>0}\subset \Lambda^3_1(M)\oplus \Lambda^3_{27}(M).
\]
Taking adjoint, we obtain also that 
\begin{equation}\label{_H^4_c_without_7-component_Equation_}
{\cal H}^4_c(M)_{>0}\subset \Lambda^4_1(M)\oplus \Lambda^4_{27}(M).
\end{equation}
Take an arbitrary $\alpha\in {\cal H}^3_c(M)_{>0}$.
Then $d\alpha\in {\cal H}^4_c(M)_{>0}$. Using
\eqref{_H^4_c_without_7-component_Equation_} and
\eqref{_eta_1_7_via_Lambda_Equation_}, we obtain that
$\Lambda d\alpha=0$. Then 
\begin{equation}\label{_d^*_c_alpha_in_Lambda_7_Equation_}
0 = d^*_c\alpha = \{\Lambda_\omega, d\} \alpha = d\Lambda_\omega\alpha.
\end{equation}
Since $\Lambda_\omega\alpha$ is a function, 
\eqref{_d^*_c_alpha_in_Lambda_7_Equation_} gives $\alpha_1=0$.
Then $\alpha\in \Lambda^3_{27}(M)$.
We proved \ref{_pseudoha_G_2_Theorem_} (ii).

Now, every $\alpha \in \Lambda^3_{27}(M)$ satisfying
$d\alpha = \mu *\alpha$ clearly belongs to 
${\cal H}^3_c(M)$. Indeed, in this case
\[ 
  L_\omega\alpha = 
C^* d\alpha = C^* \alpha = 0
\]
because the operators $C^*$, $L_\omega$ are $G_2$-invariant, and 
\[ 
  d^*\alpha = * d *\alpha = * \mu^{-1} d^2 \alpha =0.
\]
because $d^2=0$.
Taking commutators of $d^*$ with $L_\omega$
and $d^*$ with $C^*$, we find that $d_c\alpha = d^*_c\alpha=0$.
To see that such $\alpha$ generate ${\cal H}^3_c(M)$,
we use the following lemma, which finishes the proof
of \ref{_pseudoha_G_2_Theorem_} (iii).

\hfill

\lemma\label{_H^3_c_generated_by_d_alpha=*_alpha_Lemma_}
In assumptions of \ref{_pseudoha_G_2_Theorem_}, 
${\cal H}^3_c(M)_{>0}$ is generated by all
$\alpha \in {\cal H}^3_c(M)_{>0}$ which satisfy
$d\alpha = \mu *\alpha, \ \  \mu \neq 0$.

\hfill

{\bf Proof:}
Since $d_c$, $d^*_c$ commute with the Laplacian,
${\cal H}^3_c(M)_{>0}$ is generated by the eigenspaces 
${\cal H}^3_c(M)_\lambda$
of $\Delta \restrict{{\cal H}^3_c(M)_{>0}}$, which 
are finite-dimensional. Moreover, 
$* d:\; \Lambda^3(M)\arrow \Lambda^3(M)$ also commutes with the
Laplacian, hence it acts on the finite-dimensional spaces
${\cal H}^3_c(M)_\lambda$. Since
\[ 
(* d\alpha, \alpha') = \int_M d\alpha \wedge \alpha' = 
- \int_M \alpha \wedge d\alpha'= - \overline{(* d\alpha', \alpha)},
\]
the operator $* d$ is skew-Hermitian, hence semisimple. Therefore,
${\cal H}^3_c(M)_\lambda$ is generated by its eigenspaces.
By \ref{_pseudoha_G_2_Theorem_} (i),  $d^*$ vanishes
on ${\cal H}^3_c(M)$, hence $d\alpha \neq 0$ unless
$\alpha$ is harmonic. Therefore, $*d$ acts on 
${\cal H}^3_c(M)_\lambda$ with non-zero 
eigenvalues $\mu_i$.\footnote{In fact, $\lambda= |\mu_i|^2$,
as follows from $\Delta\restrict{{\cal H}^3_c(M)_{>0}} = d^* d$.}
We proved \ref{_H^3_c_generated_by_d_alpha=*_alpha_Lemma_}.
The proof of \ref{_pseudoha_G_2_Theorem_} is finished.
\endproof

\hfill

\hfill

{\bf Acknowledgements:} I am grateful to Gil Cavalcanti,
Marisa Fernandez, Nigel Hitchin,  Dima Kaledin 
and Jim Stasheff for valuable advice and 
consultations. Also, Marisa Fernandez found an error in
the original proof of \ref{_pseudoha_pseudoco_Proposition_}.

{\scriptsize

\hfill

\noindent {\sc Misha Verbitsky\\
{\sc  Institute of Theoretical and
Experimental Physics \\
B. Cheremushkinskaya, 25, Moscow, 117259, Russia }\\
\  \\
\tt verbit@maths.gla.ac.uk, \ \  verbit@mccme.ru 
}

\end{document}